\documentclass[journal]{IEEEtran}
\usepackage{amsmath,amssymb}
\usepackage{siunitx}
\usepackage{setspace,bm,comment}
\usepackage{graphicx,tabularx,multirow,float,stfloats,pdfsync}
\usepackage{graphicx,subfigure,xcolor,booktabs,makecell}
\usepackage[hidelinks]{hyperref}
\usepackage{cite,url}
\usepackage{xcolor}

\newcommand{\cM}{\mathcal{M}}
\newcommand{\E}{\mathbb{E}}
\renewcommand{\vec}[1]{\boldsymbol{#1}}
\newtheorem{theorem}{Theorem}

\begin{document}

\title{Risk-Based Distributionally Robust Optimal Power Flow With Dynamic Line Rating}

\author{Cheng~Wang, Rui~Gao, Feng~Qiu,~\IEEEmembership{ Senior Member,~IEEE}, Jianhui~Wang,~\IEEEmembership{Senior Member,~IEEE}, Linwei~Xin
\thanks{This work was supported in part by the National Natural Science Foundation of China (51725702, 51627811), and in part by the "111" project (B08013). J. Wang's work is supported by the National Science Foundation (1745451) and the U.S. Department of Energy (DOE)'s Office of Electricity Delivery and Energy Reliability.}
\thanks{C.~Wang is with the State Key Laboratory of Alternate Electrical Power System with Renewable Energy Sources, North China Electric Power University, Beijing 102206, China (e-mail: chengwang@ncepu.edu.cn).}
\thanks{R.~Gao is with H.~Milton Stewart School of Industrial \& Systems Engineering, Georgia Institute of Technology, Atlanta, GA 30332, USA (e-mail: rgao32@gatech.edu).}
\thanks{F.~Qiu is with Argonne National Laboratory, Argonne, IL 60439, USA (e-mail: fqiu@anl.gov).}
\thanks{J.~Wang is with the Department of Electrical Engineering at Southern Methodist University, Dallas, TX, USA and the Energy Systems Division at Argonne National Laboratory, Argonne, IL, USA (email: jianhui.wang@ieee.org).}
\thanks{L.~Xin is with Booth School of Business, The University of Chicago, Chicago, IL 60637, USA (e-mail: Linwei.Xin@chicagobooth.edu).}}

\maketitle

\begin{abstract}
In this paper, we propose a risk-based data-driven approach to optimal power flow (DROPF) with dynamic line rating.
The risk terms, including penalties for load shedding, wind generation curtailment and line overload, are embedded into the objective function.
To hedge against the uncertainties on wind generation data and line rating data, we consider a distributionally robust approach.
The ambiguity set is based on second-order moment and Wasserstein distance, which
 captures the correlations between wind generation outputs and line ratings, and is robust to data perturbation.
We show that the proposed DROPF model can be reformulated as a conic program.
Considering the relatively large number of constraints involved, an approximation of the proposed DROPF model is suggested, which significantly reduces the computational costs. A Wasserstein distance constrained DROPF and its tractable reformulation are also provided for practical large-scale test systems.
Simulation results on the 5-bus, the IEEE 118-bus and the Polish 2736-bus test systems validate the effectiveness of the proposed models.
\end{abstract}


\section*{Nomenclature}
Most of the symbols and notations used throughout this manuscript are defined below for quick reference. Others are defined following their first appearances, as needed.
\subsection{Sets and Indices}
\begin{IEEEdescription}[\IEEEusemathlabelsep\IEEEsetlabelwidth{$aaaaaaaaaaa$}]
\item[$g \in \mathcal{G}$]         Traditional generators.
\item[$d \in \mathcal{D}$]         Loads.
\item[$w \in \mathcal{W}$]         Wind farms.
\item[$l \in \mathcal{L}$]         Lines.
\item[$n\in\mathcal{N}$] Samples.
\item[$\mathcal{M}_1/\mathcal{M}_2/\mathcal{M}$] Ambiguity Sets.
\item[$\mathcal{P}$] The set of all probability distributions
\item[$G/W/L/N$]  Numbers of traditional generators/wind farms/lines/samples.
\end{IEEEdescription}

\subsection{Parameters}
\begin{IEEEdescription}[\IEEEusemathlabelsep\IEEEsetlabelwidth{$aaaaaaaaaa$}]
\item[$f_g(\cdot)$]                 Cost functions of traditional generators.
\item[$f_g^+(\cdot)/f_g^-(\cdot)$]  Upward/downward regulation cost functions of traditional generators.
\item[$\bar{p}_g/\underline{p}_g$] Upper/lower capacity limits of generators.
\item[$p_w$]                     Forecasting values of wind farms.
\item[$\bar{p}_w$]               Installed capacities of wind farms.
\item[$p_d$]                     Load demands.
\item[$\beta_w/\beta_d/\beta_l$]             Penalties of wind generation curtailment/load shedding/line overload.
\item[$p_l$]                     Forecasting values of lines with DLR.
\item[$\underline{p}_l$]         Static line ratings.
\item[$\pi_{gl}/\pi_{wl}/\pi_{dl}$] Shift distribution factors of generators/wind farms/loads.
\item[$\hat{\vec{\xi}}^n$] Sample vectors of wind farm outputs and line ratings.
\item[$\hat{\vec{m}}$] Mean of samples.
\item[$\hat{\vec{\Sigma}}$] Covariance of samples.
\item[$\theta$] Radius of the Wasserstein ambiguity set.
\item[$\tau$] Multiple of covariance.
\end{IEEEdescription}

\subsection{Decision Variables}
\begin{IEEEdescription}[\IEEEusemathlabelsep\IEEEsetlabelwidth{$aaaaaaaaa$}]
\item[$p_{g}$]                         Outputs of generators.
\item[$r_g^+,r_g^-$]                   Upward/downward reserves of generators.
\item[$\alpha_g$]                Reserve participation factors of generators.
\end{IEEEdescription}

\subsection{Random Variables}
\begin{IEEEdescription}[\IEEEusemathlabelsep\IEEEsetlabelwidth{$aaaaaaaaa$}]
\item[$\tilde{p}_w$]                   Actual outputs of wind farms.
\item[$\tilde{p}_l$]                   Actual line ratings of lines with DLR.
\end{IEEEdescription}

\subsection{Acronyms}
\begin{IEEEdescription}[\IEEEusemathlabelsep\IEEEsetlabelwidth{$aaaaaaaaaaa$}]
\item[AGC]                   Automatic generation control.
\item[ARO]                   Adaptive robust optimization.
\item[DC]                    Direct current.
\item[DLR]                   Dynamic line rating.
\item[DRO]                   Distributionally robust optimization.
\item[DROPF]                 Distributionally robust OPF.
\item[M-DROPF]               Moment constrained DROPF.
\item[OPF]                   Optimal power flow.
\item[QP]                    Quadratic programming.
\item[SAA]                   Sample average approximation.
\item[SDP]                   Semi-definite programming.
\item[SLR]                   Static line rating.
\item[SO]                    Stochastic optimization.
\item[UC]                    Unit commitment.
\item[W$\&$M-DROPF]          Wasserstain distance and moment constrained DROPF.
\item[W-DROPF]               Wasserstain distance constrained DROPF.
\end{IEEEdescription}

\section{Introduction}

\IEEEPARstart{W}ORLDWIDE rapid development of wind generation has brought significant economic and environmental benefits to human society, yet the concomitant uncertainty and variability also have imposed serious challenges on power system operation. One of the challenges is the lack of sufficient transmission capacity that transports wind energy from wind farms to load centers. Recently, dynamic line rating is explored to squeeze out the potential capability of existing transmission lines using new information technologies.
Traditionally, power system operation uses static line rating (SLR) (line rating indicates the maximum amount of current that the line's conductors can carry (under a set of assumed weather conditions) without violating safety codes or damaging the conductor). SLR is calculated under ``worst-case" scenarios (e.g., wind speed is 0.6 m/s, the illumination intensity is $1000 \text{W}/\text{m}^2$ and the environment temperature is $\SI{40}{\degreeCelsius}$ \cite{michiorri_forecasting_2015,fernandez_review_2016}) and tends to be very conservative. The line capacity is also sensitive to ambient environment (e.g., radiation and wind). For example, when ambient temperature decreases by $\SI{10}{\degreeCelsius}$, the line rating/capacity increases by 11\%; when wind speed perpendicular to the line increase by 1 meter per second, the capacity increase by 44\%~\cite{ValleyGroup}. Therefore, potential transmission capacity could be wasted if sticking to SLR. 
Dynamic line rating (DLR) constantly monitors the ambient environment and determines the line ratings \cite{michiorri_forecasting_2015}, which can better utilize the transmission line capacities, and is less conservative than SLR. A number of pilot projects have been implemented and the results showed that the increase of transmission capacity justified the installation costs~\cite{DemoProject}. 
Recent research on DLR includes DLR calculation and application in system operations.
\cite{carlini_reliable_2016} proposes a reliable computing framework for DLR of overhead lines. \cite{wallnerstrom_impact_2015} analyzes the impact of DLR on wind generation accommodation. \cite{jabarnejad_optimal_2016} develops a framework to determine the priority of lines to be upgraded with DLR.
The use of DLR introduces more uncertainty and might bring additional operational loss, i.e. transmission line overload \cite{ni_online_2003,li_transmission_2015}, because the actual line ratings and outputs of wind generation are random and thus cannot be accurately known, especially in look-ahead dispatch frameworks.
Besides, as both line ratings and outputs of wind generation are influenced by weather conditions, some underlying correlation may exist in the aggregated uncertainties, and such correlation is reflected by historical data.

Mathematical methods based on decision theory and robust optimization are also used to enhance the robustness of power systems against wind generation uncertainty.
For example, \cite{Ruiwei_RUC_2012} studies a two-stage robust unit commitment (UC) model using adaptive robust optimization (ARO). An effective real-time dispatch framework is established in \cite{Zhigang_AGC_2015}, in which automatic generation control (AGC) based on linear decision rule is considered. The wind generation accommodation capability of power systems under a given dispatch strategy is studied in \cite{Cheng_Assessment_2016}, which is equivalent to an ARO model mathematically. To cope with the addtional uncertainties introduced by DLR, \cite{bucher_robust_2016} provides a linear control scheme for optimal power flow (OPF) problem based on ARO theory. 

Recently, distributionally robust optimization (DRO) models are proposed, as they are data-driven and thus can explore the correlation structure of the uncertainty, and are usually less conservative than ARO.
To hedge against the line rating uncertainty, a distributionally robust congestion management model considering DLR is established in \cite{qiu_distributionally_2015}, where the impacts of uncertainties are formulated by chance constraints and the overall model is converted into a mixed integer linear programming. 
\cite{Xiong_DROUC_2017} solves a two-stage distributionally robust UC model with a moment-based ambiguity set. 
\cite{Wei_DRO_2016} and \cite{Huanhai_DRO_2016} handle the distributionally robust energy-reserve co-dispatch problem. 
\cite{zhang_distributionally_2016} and \cite{li_distributionally_2016} tackle the distributionally robust optimal power flow (DROPF) problems using a chance-constrained formulation, which guarantees the performance when the chance constraints are satisfied, but ignores the aftermath if the chance constraints are violated.

In this paper, a risk-based OPF model with DLR is proposed, using a data-driven distributionally robust approach, which hedges against data perturbations, and is able to capture the underlying correlation in the data. An ambiguity set is constructed to mathematically describe all the candidate distributions of the uncertainties, based on the empirical distribution and the statistical properties of the data. And the empirical distribution, or more precisely, the reference distribution, is formed by historical samples. Compared with existing works, the salient features of our work are summarized as below.

\begin{enumerate}
\item Our risk-based DROPF, denoted as W$\&$M-DROPF hereinafter, embeds the operational risk in the objective function, including penalties on load shedding, wind generation curtailment and line overload. To hedge against the uncertainties on wind generation data and line rating data, we consider a distributionally robust formulation in which the ambiguity set incorporates second-order moment information \cite{Xiong_DROUC_2017,Wei_DRO_2016,Huanhai_DRO_2016,li_distributionally_2016,zhang_distributionally_2016} and the Wasserstein distance to the empirical distribution \cite{gao2016distributionally,gao2017distributionally}. We are the first to consider such ambiguity set, which is suitable for hedging against high-dimensional data perturbations and can capture the correlations between wind generation outputs and line ratings.

\item A tractable reformulation of the proposed DROPF based on strong duality theory is derived, rendering a conic program. Then we suggest approximations of the reformulation, which reduce the number of constraints from exponential size to linear size in the number of AGC generators and the number of transmission lines with DLR. Such approximations have superior performance especially when the data is statistically insufficient. Then, a Wasserstein distance constrained DROPF (W-DROPF) model, and its tractable reformulation which is a convex quadratic programming (QP), are devised for practical large-scale power system from the computational tractability aspect.

\item We compare the out-of-sample performances of proposed W$\&$M-DROPF and W-DROPF models with other ones, such as sample average approximation (SAA), moment constrained DROPF (M-DROPF). Numerical results reveal the effectiveness and advantage of the proposed DROPF model over the existing ones.
\end{enumerate}

The rest of this paper is organized as follows. The mathematical formulation of the proposed DROPF model is presented in in Section II. And its tractable reformulation is given in Section III. To validate the proposed models, several numerical results on three test systems are shown in Sections IV and V. Finally, Section VI draws the conclusion.

\section{Mathematical Formulation}
In this section, major assumptions and simplifications are presented first. Then the OPF model and its distributionally robust counterpart are introduced. This section ends with the construction of the ambiguity set.

\subsection{Assumptions and Simplifications}
\begin{enumerate}
  \item Only active power related operational constraints are considered, and DC lossless power flow model is adopted \cite{roald_corrective_2016,zhang_distributionally_2016,li_distributionally_2016,bienstock_chance-constrained_2014}. 

  \item All the traditional generators are on during the dispatch periods. They all provide AGC service and their participation factors are decision variables rather than parameters. 

  \item All the transmission lines are operated with DLR mechanism. 

  \item Compared with outputs of wind generation and line ratings, the loads can be accurately known or predicted. 

  \item The proposed OPF models are single-period ones \cite{zhang_distributionally_2016,li_distributionally_2016,roald_corrective_2016,Lubin2016_chance}. Nevertheless, they can be easily extended to multi-period ones by adding subscript $t$ to each variable as well as considering the temporally-coupled ramping constraints.
\end{enumerate}

\textbf{Remark 1:} In fact, the capacity, or more accurately, the current-carring capacity, of a line is determined by its maximum designed temperature. According to \cite{IEEE_TTR}, the relationship between the accumulated heat and temperature of a bare overhead line can be described by 
\begin{equation}\label{eq:heat_steady}
q_c+q_r=q_s+I^2R(T),
\end{equation}

\noindent where (\ref{eq:heat_steady}) stands for heat balancing equations; $q_c$ and $q_r$ represent the convective and radiated heat losses, $q_s$ and $I^2R(T)$ are the solar heat gain and Joule heat terms, $T$ and $R(T)$ are the line temperature and resistance of line, respectively. It should be noted that $q_c,q_r,q_s$ can be obtained based on weather conditions, such as wind speed and solar radiation. By parameterizing $T$ in (\ref{eq:heat_steady}) with $T^{max}$, which is the maximum designed temperature, the maximum allowed ampacity of the segment of line would be known. Generally, the transmission lines are very long, and the weather conditions, which impact line ratings the most, are different along the lines. In practice, we can place weather condition monitoring devices based on certain rules, for example one device for every several towers, along the line, and the rating of each segment of the line would be known according to (\ref{eq:heat_steady}). Then we can select the minimum value as the capacity of the line, which is the line rating under DLR mechanism.

\subsection{Single-period OPF}
The mathematical formulation of the single-period OPF problem is given as below:
\begin{subequations}\label{eq:DCOPF}
\begin{align}
     \min_{p_g}\ \    &\sum_{g\in \mathcal{G}}f_g(p_g)\label{eq:obj}\\
     s.t. \ \   &\sum_{g\in \mathcal{G}}p_g+\sum_{w\in \mathcal{W}}p_w-\sum_{d\in \mathcal{D}}p_d=0,\label{eq:balance_base}\\
        & \Big|\sum_{g\in \mathcal{G}}\pi_{gl}p_g \!+ \! \sum_{w\in \mathcal{W}}\pi_{wl}p_w \!-\! \sum_{d\in \mathcal{D}}\pi_{dl}p_d\Big| \!\le\! p_l,\  \forall l\in \mathcal{\mathcal{L}},
     \label{eq:power_flow}\\
     & \underline{p}_g\le p_g \le \bar{p}_g,\ \forall g\in \mathcal{G}.\label{eq:power_capacity}
\end{align}
\end{subequations}
The objective (\ref{eq:obj}) minimizes the generation costs of the OPF problem, where $f_g(\cdot)$ is a strictly convex quadratic function. (\ref{eq:balance_base}) gives the whole network power balance condition. (\ref{eq:power_flow}) indicates the power flow limits for transmission lines. (\ref{eq:power_capacity}) presents the generation capacity limits.

\subsection{Distributionally Robust OPF}
In practice, actual outputs of wind farms and actual line ratings of transmission lines cannot be accurately known before solving problem (\ref{eq:DCOPF}), which means (\ref{eq:balance_base}) and (\ref{eq:power_flow}) may not be satisfied and undesirable operational loss might occur. In this regard, reserves are usually committed to mitigate the deviations of the random realiziations of $\tilde{p}_w,\tilde{p}_l$ from their forecast values $p_w,p_l$, which should satisfy the following constraints:
\begin{subequations}\label{eq:reserve}
\begin{align}
& 0\le r_g^+\le \bar{p}_g-p_g,\  \forall g\in \mathcal{G}\label{eq:gen_up},\\
& 0\le r_g^-\le p_g-\underline{p}_g, \ \forall g\in \mathcal{G},\label{eq:gen_down}
\end{align}
\end{subequations}
where (\ref{eq:gen_up}) and (\ref{eq:gen_down}) present boundary limits of upward and downward reserves, respectively. 
In the time scale of the OPF problem, AGC is one of the most effective reserve commitment and response mechanisms, where each generator participated in AGC affinely adjusts its output with respect to the total deviation of the uncertainties, such as the sum of forecasting error of renewables, and all the AGC generators adjust their outputs in the same direction mathematically. The mechanism of AGC can be expressed as below \cite{Jabr2013_affine}:
\begin{subequations}\label{eq:AGC}
\begin{align}
   \textstyle-r_g^{-}\le \alpha_g\sum_{w\in \mathcal{W}}(p_w-\tilde{p}_w)\le r_g^+, \  \forall g\in \mathcal{G},\label{eq:gen_capacity_AGC}\\
  0\le\alpha_g\le1,\ \ \forall g \in \mathcal{G},\label{eq:alpha_bound}\\
  \sum_{g\in \mathcal{G}}\alpha_g=1,\label{eq:alpha_sum}
\end{align}
\end{subequations}
where (\ref{eq:gen_capacity_AGC}) limits the boundaries of output adjustments of generator with respect to the deviation of sum of actual outputs of wind farms from the sum of their forecasting values, and $\alpha_g\sum_{w\in \mathcal{W}}(p_w-\tilde{p}_w)$ represents the output adjustment of generators; (\ref{eq:alpha_bound}) restricts the values of affine coefficients, a.k.a.~participation factors of AGC generators; (\ref{eq:alpha_sum}) guarantees the forecasting errors of wind generation are fully mitigated. It should be noted that the AGC generators will not adjust their outputs if the actual line rating deviates from its forecasting value in the current formulation. For line rating based generator outputs affine adjustment rules, please refer to \cite{bucher_robust_2016}. Meanwhile, the adjusted power flow should not exceed the actual line rating, resulting in for all $l\in \mathcal{L}$:
\begin{equation}\label{eq:power_flow_AGC}
  \begin{aligned}
  \big|\sum_{g\in \mathcal{G}}&\pi_{gl}\Big(p_g + \alpha_g\sum_{w\in \mathcal{W}}(p_w- \tilde{p}_w)\Big)\\ &+ \sum_{w\in \mathcal{W}}\pi_{wl}\tilde{p}_w-\sum_{d\in \mathcal{D}}\pi_{dl}p_d\big| \le \tilde{p}_l.
  \end{aligned}
\end{equation}
However, as $\tilde{p}_w$ and $\tilde{p}_l$ are random variables, which suggests (\ref{eq:gen_capacity_AGC}) and (\ref{eq:power_flow_AGC}) may not be satisfied regardless of the OPF strategy. Therefore, we penalize the expected violation of (\ref{eq:gen_capacity_AGC}) and (\ref{eq:power_flow_AGC}) and add the penalty terms into the objective function \cite{Zhang2015_risk,Wang2017_risk}, and we aim to minimize the objective under the worst-case distribution of the random variables, rendering a distributionally robust formulation of the OPF problem as below
\begin{figure*}[ht]
\begin{equation}\label{eq:DR_Obj}
\begin{split}
\min_{p_g,r_g^{+},r_g^{-},\alpha_g}\ \ \bigg\{ &\sum_{g\in \mathcal{G}}\Big(f_g(p_g)+f_g^+(r_g^+)+f_g^-(r_g^-)\Big)\\
& +\max_{\mu\in\mathcal{M}}\mathbb{E}_{\mu}\Big[\beta_d\sum_{g\in \mathcal{G}}\Big(\alpha_g\sum_{w\in \mathcal{W}}\Big(p_w-\tilde{p}_w\Big)-r_g^+\Big)^+  +\beta_w\sum_{g\in \mathcal{G}}\Big(\alpha_{g}\sum_{w\in \mathcal{W}}(\tilde{p}_w
-p_w)-r_g^-\Big)^+ \\
&\hspace{50pt} +\beta_l \sum_{l\in \mathcal{L}}\Big(\big|\sum_{g\in \mathcal{G}}\pi_{gl}\Big(p_g+\alpha_g\sum_{w\in \mathcal{W}}\Big(p_w-\tilde{p}_w\Big)\Big)+\sum_{w\in W}\pi_{wl}\tilde{p}_w-\sum_{d\in \mathcal{D}}\pi_{dl}p_d\big|-\tilde{p}_{l}\Big)^+\Big]\bigg\}
\end{split}\tag{6}
\end{equation}
\hrulefill
\end{figure*}
\setcounter{equation}{6}
\begin{equation}\label{problem:DROPF}
  \begin{aligned}
    \textrm{Objective: } &   (\ref{eq:DR_Obj})  \\
    \textrm{Constraints: } &   \text{(\ref{eq:balance_base})-(\ref{eq:power_capacity}),\ (\ref{eq:reserve}),\ (\ref{eq:alpha_bound})-(\ref{eq:alpha_sum})},
  \end{aligned} \tag{\textsf{W$\&$M-DROPF}}
\end{equation}
where (\ref{eq:DR_Obj}) serves as the objective function of the proposed model, the first term is the same as (\ref{eq:obj}); the second and third terms are the regulation costs for committing upward and downward reserves, respectively, where both $f_g^+(\cdot)$ and $f_g^-(\cdot)$ are strictly convex quadratic functions; the rest terms are the penalties for load shedding, wind generation curtailment and line overload, respectively; and $\mu$ is any distribution in the ambiguity set $\mathcal{M}$; $\mathbb{E}$ is the expectation operator; $(\cdot)^+$ returns the larger value between zero and the whole expression in the brackets; $|\cdot|$ takes the absolute value of the internal term. In (\ref{problem:DROPF}), constraints (\ref{eq:balance_base})-(\ref{eq:power_capacity}) ensure the existence of a feasible OPF strategy when the outputs of wind generation and line ratings take their predicted values, where neither upward nor downward reserve is committed; constraints (\ref{eq:reserve}) and (\ref{eq:alpha_bound})-(\ref{eq:alpha_sum}) render practical reserve allocation decisions. It should be noted that power imbalance and line overload issues may occur if the OPF strategy obtained from the constraints of (\ref{problem:DROPF}) is adopted. Specifically, the power imbalance issue originates the violation of (\ref{eq:gen_capacity_AGC}), where the violations of the left- and right-side inequality lead to load shedding and wind generation curtailment, respectively; and the line overload issue is on account of the violation of (\ref{eq:power_flow_AGC}). Therefore, the risk terms representing load shedding, wind generation curtailment as well as line overload, are added in the objective function of (\ref{problem:DROPF}) to improve the performance of the OPF strategy.

In (\ref{problem:DROPF}), the operation risk terms are optimized along with the operation costs in the objective function, however, it can be easily extended to a risk-limiting form \cite{Risk_limiting,CC_Weighted} by removing the risk related terms in the objective function and adding a risk limit constraint instead.

\textbf{Remark 2:} There are many inspiring works on dealing with the uncertainties in power system operation problems \cite{Stochastic_Jianhui,Stochastic_Qianfan,Stochastic_Wu,bienstock_chance-constrained_2014,Ruiwei_RUC_2012,Cheng_Assessment_2016,Zhigang_AGC_2015,Xiong_DROUC_2017,Wei_DRO_2016,Huanhai_DRO_2016}, where the two-stage modelling framework is commonly adopted. And they can be grouped into three categories according to the mathematical formulations, which are

1) Stochastic optimization (SO) based works \cite{Stochastic_Jianhui,Stochastic_Qianfan,Stochastic_Wu,bienstock_chance-constrained_2014}, which can be expressed by the general form
\begin{equation}\label{eq:stochastic_com}
O_S=\min_{\vec{x}\in\mathcal{X}}\vec{c}^\top\vec{x}+\mathbb{E}_{\vec{\xi}\in\mu}\left[\min_{\vec{y}\in\mathcal{Y}(\vec{x},\vec{\xi})}\left(\vec{d}^\top\vec{\xi}+\vec{e}^\top\vec{y}\right)\right],
\end{equation}
where $\vec{x}$ and $\vec{y}$ are the first- and second-stage decision variables, respectively, or referred to as \textit{here-and-now} and \textit{wait-and-see} decisions, respectively; $\vec{\xi}$ is the realization of the uncertainty, and $\mu$ is its distribution, which is assumed to be known beforehand; $\vec{c},\vec{d},\vec{e}$ are constant coefficient vectors; $\mathcal{X}$ and $\mathcal{Y}$ are the feasible regions for $\vec{x}$ and $\vec{y}$, respectively; $\mathbb{E}$ is the expectation operator; $O_S$ is the value of the (\ref{eq:stochastic_com}), which is the sum of the first-stage decision costs and the expectation of the second-stage decision costs.

2) ARO based works \cite{Ruiwei_RUC_2012,Cheng_Assessment_2016,Zhigang_AGC_2015}, whose general expression is given as below 
\begin{equation}\label{eq:robust_com}
O_R=\min_{\vec{x}\in\mathcal{X}}\vec{c}^\top\vec{x}+\max_{\vec{\xi}\in\mathcal{U}}\min_{\vec{y}\in\mathcal{Y}(\vec{x},\vec{\xi})}\left(\vec{d}^\top\vec{\xi}+\vec{e}^\top\vec{y}\right),
\end{equation}
where $\mathcal{U}$ is a pre-determined uncertainty set, and definitions of the rest terms are the same with (\ref{eq:stochastic_com}). $O_R$ denotes the value of (\ref{eq:robust_com}), which is the sum of the first-stage decision costs and the worst-case second-stage decision costs.

3) DRO based works \cite{Xiong_DROUC_2017,Wei_DRO_2016,Huanhai_DRO_2016}, who take the form of
\begin{equation}\label{eq:distributionally_com}
O_{DR}=\min_{\vec{x}\in\mathcal{X}}\vec{c}^\top\vec{x}+\max_{\mu\in\mathcal{M}}\mathbb{E}_{\vec{\xi}\in\mu}\left[\min_{\vec{y}\in\mathcal{Y}(\vec{x},\vec{\xi})}\left(\vec{d}^\top\vec{\xi}+\vec{e}^\top\vec{y}\right)\right],
\end{equation}
where $\mathcal{M}$ is a pre-determined distribution set and $\mu$ is the worst-case distribution. We use $O_{DR}$ to represent the value of (\ref{eq:distributionally_com}), which is the sum of the first-stage decision costs and the expectation of the second-stage decision costs under the worst-case distribution. $\mathcal{M}$ is usually referred to as the ambiguity set, and obviously it has a direct influence on $O_{DR}$ and the out-of-sample performance of $\vec{x}$.

\subsection{Ambiguity Set Construction}
Before tackling the proposed model (\ref{problem:DROPF}), it is crucial to construct a meaningful and tractable ambiguity set. Let $\nu :=\frac{1}{N}\sum_{n=1}^{N}\delta_{\hat{\vec{\xi}}^n}$ be the empirical distribution, where $\hat{\vec{\xi}}^n,n=1,\ldots,N$ are samples and $\delta_{\hat{\vec{\xi}}^n}$ represents the Dirac measure on $\hat{\vec{\xi}}^n$. Two main aspects are mainly considered during the construction of $\mathcal{M}$: (i) any distribution $\mu\in\mathcal{M}$ should be close to the empirical distribution $\nu$ in the sense of proper statistical distance; (ii) any distribution $\mu\in\mathcal{M}$ should also have a similar correlation structure as $\nu$ does.

To capture desiderata (i), given $\nu$, we define
\[\label{eq:wasserstein_set}
  \mathcal{M}_1:=\Big\{\mu\in\mathcal{P}:\Theta(\mu,\nu)\le\theta\Big\},
\]
where $\theta$ is a positive parameter, and $\Theta(\mu,\nu)$ is the Wasserstein distance (of order 1) between $\mu$ and $\nu$, given by
\[\label{eq:wasserstein_def}
  W(\mu,\nu):=\min_{\gamma}\Big\{\int_{\mathbb{R}^{(W+L)}\times\mathbb{R}^{(W+L)}}\parallel\vec{\xi}-\vec{\zeta}\parallel\gamma\Big(\text{d}\vec{\xi},\text{d}\vec{\zeta}\Big)\Big\},
\]
where $\gamma$ is a joint distribution on $\mathbb{R}^{(W+L)}\times\mathbb{R}^{(W+L)}$ with marginals $\mu,\nu$; $\vec{\xi},\vec{\zeta}$ are the integral variables; $\parallel\cdot\parallel$ is the norm operator. Thus $\mathcal{M}_1$ contains all probability distributions whose Wasserstein distancees to the empirical distribution are no more than $\theta$. Wasserstein distance is well suited for hedging against the perturbation of data values and has good out-of-sample performance \cite{esfahani2015data,gao2016distributionally}.
This is particular useful to study large-scale OPF problems which involve a great number of generators and tranmission lines.

To capture desiderata (ii), we define
\[
  \mathcal{M}_2:=\Big\{\mu\in\mathcal{P}:\mathbb{E}_{\mu}[(\vec{\xi}-\hat{\vec{m}})(\vec{\xi}-\hat{\vec{m}})^\top]\preceq\tau\hat{\vec{\Sigma}}\Big\},
\]
where $\tau\ge1$, $\hat{\vec{m}}$ and $\hat{\vec{\Sigma}}$ are the mean and covariance matrix of the empirical distribution $\nu$ generated by the samples $\hat{\vec{\xi}}^n$; $\vec{\xi}$ expresses the random variables. The constraint basically suggests that the centered second-moment matrix of any relevant distribution, which reflects the correlation structure, should be close to that of the empirical distribution.
Finally, we set
\begin{equation}\label{eq:uncertainty_set}
  \mathcal{M}:=\mathcal{M}_1\cap\mathcal{M}_2.
\end{equation}
The ambiguity set $\cM$ contains all probability distributions that are close to the empirical distribution, and has a similar correlation structure to the empirical distribution.
Since $\mathcal{M}$ is comprised of infinitely many distributions, the proposed \eqref{problem:DROPF} model is not immediately computationally tractable, and we will provide a tractable formulation in the next section.

\textbf{Remark 3:} As mentioned in the Introduction, the empirical distribution $\nu$ is formed by historical samples. It should be noted that the candidate sample set for $\nu$ only consists a relatively small part of all the historical data, whose meteorological conditions are similar with the ones of the current decision-making stage. Considering the relatively high dimensionality of the random variable vector, which consists the outputs of wind farms and the line ratings, there may not be too many available samples in practical applications. Some existing works have similar settings for the number of available samples. In \cite{zhang_distributionally_2016}, they assume the decision-maker has limited knowledge of the uncertainties and set the sample number $N=20$.

\section{Solution Methodology}
In this section, we first derive a tractable reformulation of the proposed \eqref{problem:DROPF} model using strong duality theory. Then an approximation of the reformulation is suggested in order to speed up the computation. At last, the W-DROPF and its tractable reformulation are proposed for practical large-scale power systems, which reflects the tradeoff between the out-of-sample performance and the computational burden of the model.

\subsection{Conic Program Reformulation}
To simplify notations in the model \eqref{problem:DROPF}, we rewrite it as
\begin{subequations}\label{eq:DRcom}
\begin{align}
  \min_{\vec{x}}\ &  \vec{f}(\vec{x})+\max_{\mu\in\mathcal{M}}\mathbb{E}_{\mu}\left[\Psi(\vec{x},\vec{\xi})\right]\label{eq:DRcom_Obj}\\
  s.t.\  & \vec{Ax}\le\vec{h}\label{eq:DRcom_Con}
\end{align}
\end{subequations}
\noindent where $\vec{x}=[(p_g;r_g^+;r_g^-;\alpha_g): g\in \mathcal{G}]$ is a vector of decision variables; $\vec{\xi}=[(\tilde{p}_w:w\in \mathcal{W});(\tilde{p}_l:l\in \mathcal{L})]$ is a vector of random variables; $\vec{A},\vec{h}$ are coefficient matrices and can be derived from the constraints (\ref{eq:DCOPF}), (\ref{eq:reserve}) and (\ref{eq:alpha_bound})-(\ref{eq:alpha_sum}) of \eqref{problem:DROPF}; $\vec{f}(\cdot)$ expresses the first summation in \eqref{eq:DR_Obj}; and 
\begin{equation}\label{eq:Risk_term}
\max_{\mu\in\mathcal{M}}\mathbb{E}_{\mu}\left[\Psi(\vec{x},\vec{\xi})\right]
\end{equation}
represents the penalties under the worst-case distribution; and
\begin{equation}\label{eq:Psi}
  \Psi(\vec{x},\vec{\xi}) = \max_{1\leq k\leq K} \vec{a}_k(\vec{x})^\top \vec{\xi} + b_k(\vec{x})
\end{equation}
expresses the piecewise-linear convex function of $\vec{\xi}$ inside the expectation $\mathbb{E}_\mu$ in \eqref{eq:DR_Obj}, where $k$ is the index for the piece-wise linear segment of (\ref{eq:Psi}), $K=4^G\times3^L$, $\vec{a}_k(\vec{x}),b_k(\vec{x})$ are the coefficients that can be derived from the objective function \eqref{eq:DR_Obj} and their detailed expressions are provided in the Appendix.A.

Using Corollary 1 in \cite{gao2017distributionally}, 
we obtain the following theorem.

\begin{theorem}\label{thm:dual}
  Problem (\ref{eq:DRcom}) admits a conic program reformulation \eqref{eq:strong_dual}.
  \begin{figure*}[t]
    \begin{subequations}\label{eq:strong_dual}
      \begin{align}
          \min_{\substack{\vec{x},\vec{y},\vec{\zeta}\\ \lambda\geq0,\vec{\Gamma}\succeq0}} \quad & \vec{f}(\vec{x}) + \lambda\theta + \tau\text{tr}(\vec{\Gamma}\hat{\vec{\Sigma}}) + \frac{1}{N}\sum_{n=1}^N y_n\label{eq:dual_obj}\\
          s.t.\quad & \left[\begin{array}{cc}
            \vec{\Gamma} & -\vec{a}_k(\vec{x})/2+\vec{\zeta}^{nk}/2-\vec{\Gamma}\hat{\vec{m}}\\
            (-\vec{a}_k(\vec{x})/2+\vec{\zeta}^{nk}/2-\vec{\Gamma}\hat{\vec{m}})^\top\ \  & y_n-b_k(\vec{x})-(\vec{\zeta}^{nk})^\top\hat{\vec{\xi}}^n+\hat{\vec{m}}^\top\vec{\Gamma}\hat{\vec{m}}
          \end{array}\right] \succeq 0,\ \ \forall 1\leq k\leq K,\ n\in\mathcal{N},\label{eq:dual_SDP}\\
          & ||\vec{\zeta}^{nk}||_\ast \leq \lambda,\ \ \forall 1\leq k\leq K,\ n\in\mathcal{N},\label{eq:dual_SOCP}\\
          & \vec{Ax}\le\vec{h}.\label{eq:DRcom_Con_2}
      \end{align}
    \end{subequations}
    \hrule
  \end{figure*}
\end{theorem}

In reformulation \eqref{eq:strong_dual}, $y_n,n=1,...N$ and $\lambda$ are auxiliary variables; $\vec{\zeta}^{nk}\in\mathbb{R}^{(L+W)}, n=1,...,N, k=1,...,K$ are auxiliary variable vectors; $\vec{\Gamma}\in\mathbb{R}^{(L+W)}\times\mathbb{R}^{(L+W)}$ is an auxiliary variable matrix; tr is the trace operator; $||\cdot||_\ast$ is the norm dual to the norm in the definition of Wasserstein distance; constraints \eqref{eq:dual_SDP} are semi-definite programming (SDP) constraints; constraints (\ref{eq:dual_SOCP}) and (\ref{eq:DRcom_Con_2}) are linear constraints, and (\ref{eq:DRcom_Con_2}) is identical to (\ref{eq:DRcom_Con}). The objective function (\ref{eq:dual_obj}) is convex and quadratic. Therefore, (\ref{eq:strong_dual}) suggests an SDP and can readily be solved by the off-the-shelf solvers such as MOSEK. 

For large instances, (\ref{problem:DROPF}) might be computationally challenging, as it has exponentially many ($4^G\times3^L\times N+1$) SDP constraints. And SDP constraints reduction is a straightforward method to decrease the computational burden. One can adopt the inactive transmission line capacity constraints identification method in \cite{Hongxing_Line_2016} to reduce the value of $L$, where the line rating should be $\underline{p}_l$ rather than $p_l$, guaranteeing the identified inactive line capacity constraints would always stay inactive in practice regardless of the actual line ratings, and then the number of SDP constraints would become smaller. An underlying assumption is that the SLR $\underline{p}_l$ is the lower bound of $p_l$, which can be supported by the fact that SLR is usually calculated under ``worst-case" weather conditions and tends to be very conservative. Besides, the approximation of penalty term (\ref{eq:Risk_term}) would also bring similar computational benefits by reducing the total number of SDP constraints, which would be introduced in the next subsection.

\subsection{Penalty Terms Approximation}
We write $\Psi(\vec{x},\vec{\xi})$ as the sum of three terms
$\Psi(\vec{x},\vec{\xi}) = \Psi_1(\vec{x},\vec{\xi})+ \Psi_2(\vec{x},\vec{\xi})+ \Psi_3(\vec{x},\vec{\xi})$,
where the three terms on the right-hand side of the equation above represent the three summations in \eqref{eq:DR_Obj}, representing the risk of load shedding, wind generation curtailment and line overload, respectively.
By interchanging the maximum and summation, the following quantity
\begin{equation}\label{eq:DROPF_approximate_small}
  \max_{\mu\in\cM}\;\E_{\mu}[\Psi_1(\vec{x},\vec{\xi})] + \max_{\mu\in\cM}\;\E_{\mu}[\Psi_2(\vec{x},\vec{\xi})] + \max_{\mu\in\cM}\;\E_{\mu}[\Psi_3(\vec{x},\vec{\xi})]
\end{equation}
provides an upper bound on (\ref{eq:Risk_term}).
The reformulation of the corresponding \eqref{problem:DROPF} problem can be obtained similarly to Theorem \ref{thm:dual}.
After such approximation, the computational costs of the proposed \eqref{problem:DROPF} model can be reduced, as the number of SDP constraints decreases to $\left((2^{(G+1)}+3^L)\times N+3\right)$.

We can further reduce the size of the problem by considering a cruder approximation.
We formulate distributionally robust problems for each individual generator and each transmission line with DLR separately.
More specifically, consider
\begin{equation}\label{eq:DROPF_approximate_big}
  \begin{aligned}
  \sum_{g \in \mathcal{G}} & \Big(\max_{\mu\in\cM}\E_{\mu}[\Psi_1^g(\vec{x},\vec{\xi})]+\max_{\mu\in\cM}\E_{\mu}[\Psi_2^g(\vec{x},\vec{\xi})]\Big)\\
  &+\sum_{l\in \mathcal{L}}\max_{\mu\in\cM}\E_{\mu}[\Psi_3^l(\vec{x},\vec{\xi})],
  \end{aligned}
\end{equation}
where $\Psi_1^g,\Psi_2^g,\Psi_3^l$ are the summands inside the inner maximization problem in \eqref{eq:DR_Obj}.
Such approximation yields a reformulation with $\left((4G+3L)\times N+2G+L\right)$ number of SDP constraints, linearly growing with $G$ and $L$.
In fact, (\ref{eq:DROPF_approximate_big}) suggests a modelling choice, where the worst-case distributions are specified for each individual random variable related constraint, rather than one worst-case distribution for all the random variable related constraints. It should be noted that similar modelling choice can be found in the literature. Most of them are chance-constrained optimization problems, e.g., in \cite{zhang_distributionally_2016,Lubin2016_chance,bienstock_chance-constrained_2014}.

We note that (\ref{eq:DROPF_approximate_small}) and (\ref{eq:DROPF_approximate_big}) provide upper bounds on the operational risk in (\ref{eq:Risk_term}), which means they would bring additional conservativeness to the original DROPF problem. When the number of available samples is relatively small, the parameters in the ambiguity set $\mathcal{M}$, say the Wasserstein radius $\theta$ and the covariance multiple $\tau$, may not be properly tuned, which may worsen the out-of-sample performance of the proposed model. In this situation, the additional conservativeness introduced by penalty term approximation might be beneficial to the out-of-sample performance of the proposed model, as will be demonstrate by the simulation results in Section IV.D.

\subsection{W-DROPF and Its Tractable Reformulation}
For practical large-scale power systems, the proposed W$\&$M-DROPF model may not be applied even if the risk terms in its objective function are replaced with the approximated expression (\ref{eq:DROPF_approximate_big}), as the number of SDPs is still large nad the SDP constraints usually require large computational efforts. To cope with the computability issue, a DROPF model with Wasserstein distance constraint in the ambiguity set, denoted as W-DROPF, is developed and its tractable reformulation falls into a convex QP, indicating a acceptable computational time for large-scale test systems.

It should be noted that the only difference between the W$\&$M-DROPF and W-DROPF models lies in the ambiguity set, where the ambiguity set of the former model consists both Wasserstein distance and second-order moment constraints, say $\mathcal{M}=\mathcal{M}_1\cap\mathcal{M}_2$, and the one of the latter model only consists Wasserstein distance constraint, say $\mathcal{M}_1$. Therefore, the detailed W-DROPF model is not listed for simplicity. And according to Corollary 2 in \cite{gao2017distributionally}, W-DROPF has a strong dual reformulation as below
\begin{subequations}\label{eq:DROPF_only_wasserstein}
\begin{align}
&\min_{\substack{\vec{x},\vec{y}\\ \lambda\ge 0}} \ \vec{f}(\vec{x})+\lambda\theta+\frac{1}{N}\sum_{n=1}^{N}y_n\label{eq:DROPF_only_wasserstein_obj}\\
s.t.\ \ & y_n\ge \vec{a}_k(\vec{x})^\top\vec{\hat{\xi}}^n+b_k(\vec{x}),\ \forall 1\le k\le K, \ n\in\mathcal{N},\label{eq:DROPF_only_wasserstein_con1}\\
\begin{split}
&\lambda\ge\max\left\{||\vec{a}_{j_{g_d}}^{g_d}(\vec{x})||_\ast,\ ||\vec{a}_{j_{g_w}}^{g_w}(\vec{x})||_\ast,\ ||\vec{a}_{j_l}^l(\vec{x})||_\ast \right\},\\
& \forall g\in \mathcal{G},\ l\in\mathcal{L},\ j_{g_d},j_{g_w}\in\{1,2\},\ j_l\in\{1,2,3\},
\end{split}\label{eq:DROPF_only_wasserstein_con2}\\
& \vec{Ax}\le\vec{h},\label{eq:DROPF_only_wasserstein_con3}
\end{align}
\end{subequations}

\noindent where $\vec{x}=[(p_g;r_g^+;r_g^-;\alpha_g): g\in \mathcal{G}]$ is a vector of decision variables; $y_n,n=1,...N$ and $\lambda$ are auxiliary variables; $\vec{a}_k(\vec{x}),b_k(\vec{x}),\vec{a}_{j_{g_d}}^{g_d}(\vec{x}),\vec{a}_{j_{g_w}}^{g_w}(\vec{x}),\vec{a}_{j_l}^l(\vec{x})$ are the coefficients and their expressions can be found in Appendix.A; $||\cdot||_\ast$ is the norm dual to the norm in the definition of Wasserstein distance. In reformulation (\ref{eq:DROPF_only_wasserstein}), the objective function (\ref{eq:DROPF_only_wasserstein_obj}) is quadratic and convex, and the constraints (\ref{eq:DROPF_only_wasserstein_con1})-(\ref{eq:DROPF_only_wasserstein_con3}) are linear constraints, and (\ref{eq:DROPF_only_wasserstein_con3}) is identical to (\ref{eq:DRcom_Con}). Evidently, (\ref{eq:DROPF_only_wasserstein}) suggests a convex QP, which means it can be efficiently solved by commercial solvers such as Cplex and Gurobi.

\section{Illustrative Example}
In this section, we present numerical results on a 5-bus test system to validate the proposed methods. All experiments are performed on a laptop with Intel\textregistered~Core\texttrademark~2 Duo 2.2 GHz CPU and 4 GB memory. The proposed algorithms are coded in MATLAB with YALMIP toolbox \cite{Lofberg2004}. SDPs are solved by MOSEK, while QPs are solved with Gurobi. Four different models are tested and we list them as below for a quick reference: 1) W$\&$M-DROPF, denoted as A1. 2) W-DROPF, denoted as A2. 3) M-DROPF, denoted as A3. 4) SAA, denoted as A4, where the tractable reformulations of A3 suggests an SDP and is presented in the Appendix.A, and A4 is a QP. It should be noted A1-A4 share the same primal decision variables ($p_g,r^+_g,r^-_g,\alpha_g$) and constraints ((\ref{eq:balance_base})-(\ref{eq:power_capacity}), (\ref{eq:gen_up})-(\ref{eq:gen_down}), and (\ref{eq:alpha_bound})-(\ref{eq:alpha_sum})). And the numbers of auxiliary variables and constraints generated during model reformulation and approximation are summarized in Table \ref{Tab:computational_scale_general}, where $N,G,L,W$ represent the numbers of samples, traditional generators, lines, and wind farms, respectively. As the objective approximation methods proposed in Section III.B can also be applied in A2 and A3, the corresponding numbers of auxiliary variables and constraints are also listed in Table \ref{Tab:computational_scale_general}. From Table \ref{Tab:computational_scale_general}, it can be observed that A2 and A4 are still QPs after model reformulation and approximation, while A1 and A3 become SDPs.
\begin{table*}[ht!]
\footnotesize
  \centering
  \caption{Number of auxiliary constraints and variables in the proposed models.}\label{Tab:computational_scale_general}
  \begin{tabular}{ccccc}
   \toprule
   \multicolumn{2}{c}{\multirow{2}*{Model Type}} & \multirow{2}*{Auxiliary Variable} & \multicolumn{2}{c}{Auxiliary Constraint}  \\ 
   \cmidrule(lr){4-5}
   & & & SDP & Linear\\
   \midrule
   \multirow{3}*{A1: W$\&$M-DROPF} & (\ref{eq:Risk_term})  & $(L+W)\times 	N\times 4^G\times 3^L+N+2$ & $N\times 4^G\times 3^G+1$ & $N\times 4^G\times 3^G+1$\\
   & (\ref{eq:DROPF_approximate_small})  & $(L+W)\left(2^{G+1}+3^L\right)\times N+3N+6$ & $(2^{G+1}+3^L)N+3$ & $(2^{G+1}+3^L)N+3$\\
   & (\ref{eq:DROPF_approximate_big}) & $(N+2)(2G+L)+(L+W)(4G+3L)N$ & $(4G+3L)N+2G+L$ & $(4G+3L)N+2G+L$ \\
   \midrule
   \multirow{3}*{A2: W-DROPF}& (\ref{eq:Risk_term}) & $N+1$ & $0$ & $N\times 4^G\times 3^L+2G+2L+1$ \\
   & (\ref{eq:DROPF_approximate_small}) & $3N+3$ & $0$ & $(2^{G+1}+3^L)N+6G+6L$\\
   & (\ref{eq:DROPF_approximate_big}) & $(2G+L)(N+1)$ & $0$ & $(4G+3L)N+2G+2L$\\
   \midrule
   \multirow{3}*{A3: M-DROPF}& (\ref{eq:Risk_term}) &$(L+W)^2+L+W+1$& $4^G\times3^L+1$&$0$  \\
   & (\ref{eq:DROPF_approximate_small}) & $3((L+W)^2+L+W+1)$ & $2^{G+1}+3^L+3$& $0$\\ 
   & (\ref{eq:DROPF_approximate_big}) & $((L+W)^2+L+W+1)(2G+L)$ & $6G+4L$& $0$\\
   \midrule
   \multicolumn{2}{c}{A4: SAA} & $2NG+NL$& $0$ & $4NG+3NL$ \\ 
   \bottomrule
  \end{tabular}
\end{table*}

\subsection{Test System Description}
\begin{figure}[ht]
\centering
  \includegraphics[width=0.4\textwidth]{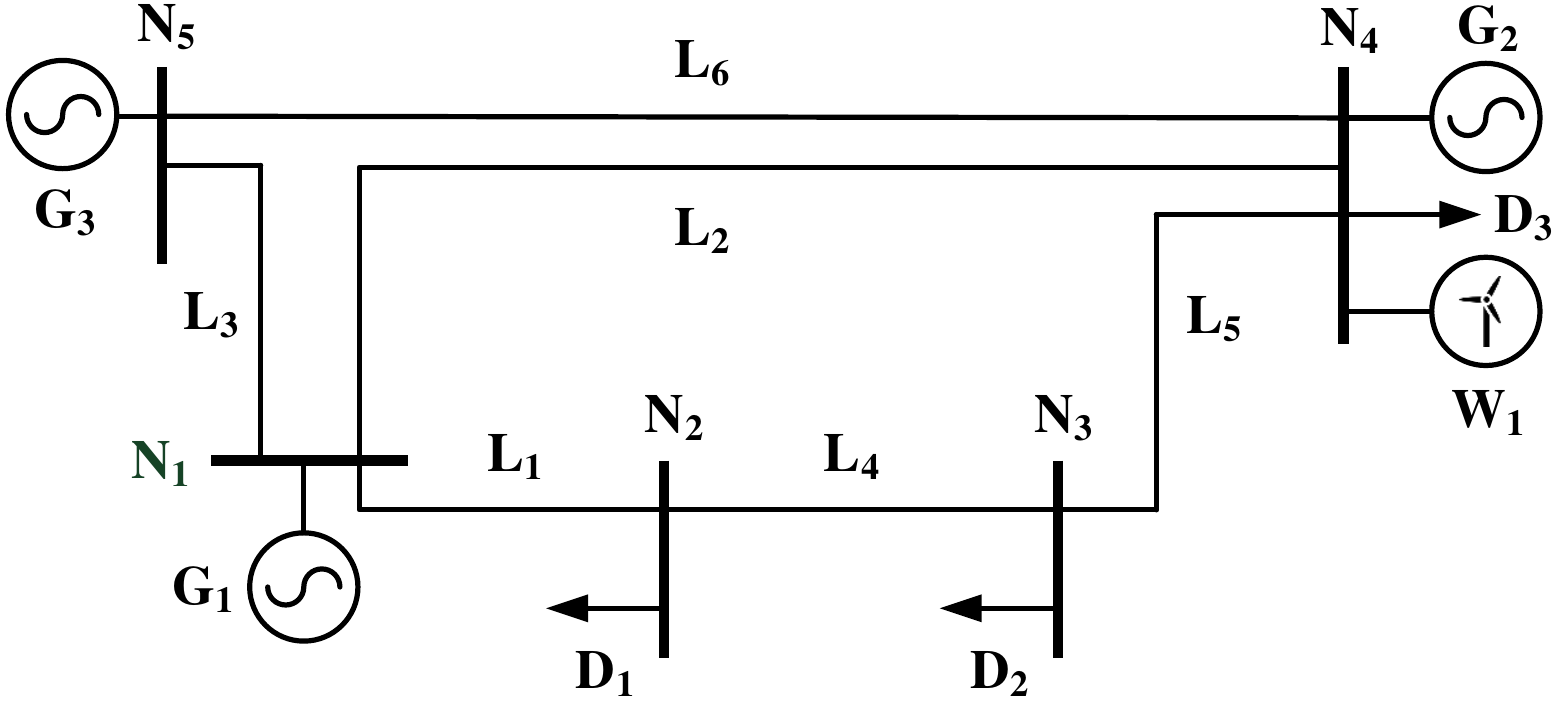}
  \caption{Topology of the 5-bus test system.}
  \label{fig:Topology_5bus}
\end{figure}

Fig. \ref{fig:Topology_5bus} depicts the topology of the test system, which will be referred to as the 5-bus test system later on. The 5-bus test system has 5 buses, 3 traditional generators, 1 wind farm, 6 transmission lines and 3 loads, which are denoted by N, G, W, L and D with subscripts, respectively. All the traditional generators are AGC units, and are always on during dispatch periods. Please refer to \cite{5bus} for the load profile, the penalty coefficients and detailed system data. The forecasting value of line with DLR, i.e., $p_l$, is $20\%$ larger than the SLR $\underline{p}_l$ in all the cases. Besides, the line ratings listed in \cite{5bus} are SLRs.

\subsection{Sample Generation \& Parameter Tuning}
The inactive transmission line capacity constraints identification method in \cite{Hongxing_Line_2016} is adopted, and $\text{L}_2,\text{L}_3,\text{L}_4$ are recognized as redundant lines and their corresponding capacity limits are neglected in the following analysis. Thus, the remaining random variables in the proposed DROPF model are the outputs of $\text{W}_1$ and line ratings of $\text{L}_1,\text{L}_5,\text{L}_6$. And the dimensionality of the uncertainty vector $\vec{\xi}=[(\tilde{p}_w:w\in\mathcal{W});(\tilde{p}_l:l\in\mathcal{L})]$ of the 5-bus test system is 4.

We assume $\vec{\xi}$ follows a multivariate Gaussian distribution, where the mean values of $\vec{\xi}$ can be found in \cite{5bus} and the standard deviation of each element of $\vec{\xi}$ equals its mean value multiplies a parameter randomly chosen from $[0.5,1]$. To guarantee the validity of the samples, a sample validity checking procedure is developed and executed during sample generation. Specifically, for outputs of wind farms, valid samples should satisfy
\begin{equation}\label{eq:valid_wind}
0\le\tilde{p}_w^n\le\bar{p}_w, \ \forall w\in\mathcal{W},\ \forall n\in\mathcal{N},
\end{equation}
where $\tilde{p}_w^n$ is the outputs of wind farm $w$ in sample $\hat{\vec{\xi}}^n$; $\bar{p}_w$ is the installed capacity of wind farm $w$. Similarly, for line ratings, valid samples should meet
\begin{equation}\label{eq:valid_line}
\underline{p}_l\le\tilde{p}_l^n, \ \forall l\in\mathcal{L},\ \forall n\in\mathcal{N},
\end{equation} 
where $\tilde{p}_l^n$ is the rating of line $l$ in sample $\hat{\vec{\xi}}^n$; $\underline{p}_l$ is the static rating of line $l$. As mentioned in the introduction, SLR is calculated under ``worst-case" weather conditions and tends to be very conservative. Thus, the value of SLR should be the lower bound of line capacity. It should be noted that the proposed sample validation procedure is always executed during sample generation unless specified.

The correlation coefficient between $\xi_i$ and $\xi_j$ are modeled as $\rho^{|i-j|}$ \cite{zou2005regularization}, where $\rho$ is chosen from $\{0.4,0.9\}$ and $i,j$ are the indices for the elements in $\vec{\xi}$. Note that $0.4^2= 0.16$ and $0.4^3\simeq 0.06$, which means most of the components have small correlations. Similary, $0.9^2=0.81$ and $0.9^3\simeq 0.73$, which means most of the components have large correlations. We assume the knowledge of the decision maker over the uncertainty $\vec{\xi}$ is limited, and a small sample set with $N=20$ is considered in the following analysis. The simulation is repeated by 50 times for a specific $\rho$.

The tuning parameters, i.e., Wasserstein radius $\theta$, is selected using hold-out cross validation method. In each repetition, the $N$ samples are randomly partitioned into a training dataset (70\% of the samples) and a validation dataset (the remaining 30\%). For different tuning parameters, we use the training dataset to solve problem (\ref{eq:strong_dual}) and then use the validation dataset to estimate the out-of-sample performance of different parameter values and select the one with the best performance, i.e. the optimal objective value in (\ref{eq:strong_dual}). The we resolve problem (\ref{eq:strong_dual}) with the best tuning parameters using all $N$ samples and obtain the optimal solution. Finally, the performances of A1-A4 are examined using an independent testing dataset consists of $10^4$ samples. Specially, the ranges for Wasserstein radius $\theta$ and covariance matrix multiple $\tau$ tuning are $[0,\ 1]$ and $[1,\ 10]$, respectively. The tuning steps for $\theta$ and $\tau$ are $0.01$ and $1$, respectively. From simulation results, $\theta$ and $\tau$ take different values in different sample sets. The mean values of $\theta$ and $\tau$ are $0.0527$ and $2.09$, respectively, for the 100 sample sets used in Section IV.C.

\subsection{Simulation Results}
In the sequel, detailed comparisons are made among the aforementioned approaches, i.e., A1-A4, for the data-driven OPF problem. Besides, to verify the significance of considering the DLR mechanism for transmission lines, we add control groups for A1-A4, where all the lines are operated with SLRs and the rest of the model for a specific approach is the same. Table \ref{Tab:comparison_value_5bus} collects all the results in the field of effectiveness, where DC and OP are short for dispatch costs (sum of first three terms in (\ref{eq:DR_Obj})) and out-of-sample performances (dispatch costs plus operational risk under the testing dataset), respectively. Besides, in Table \ref{Tab:comparison_value_5bus}, avg, max and min denote the average, worst and best performances of an approach, respectively.

The effectiveness of DLR mechanism is first investigated. It can be observed that all the values of the cost terms of OPF strategies, i.e., DC in Table \ref{Tab:comparison_value_5bus}, of A1-A4 under different $\rho$ with DLR mechanism are lower than those with SLR mechanism. The reason is that SLRs are always conservative and adopting DLR mechanism can enlarge the feasible region of the OPF problem, resulting in the decrement in dispatch costs. However, the gaps of average out-of-sample performance between DLR and SLR mechanisms get smaller compared with the corresponding dispatch costs gaps. Note that the worst performance of A1 when $\rho=0.9$ with DLR mechanism is worse than that with SLR mechanism, and similar observation can be found in the best performance of A1 when $\rho=0.4$. This is because the transmission line overload risk of under SLR mechanism might be smaller than that under DLR mechanism. Nevertheless, the average performances with DLR mechanism are always better than those with SLR mechanism, which indicates the effectiveness of the DLR mechanism and the data-driven OPF models can effectively hedge against the uncertainties of line ratings.

Then we will focus on the comparisons among the uncertainty modeling approaches, i.e. A1-A4. From Table \ref{Tab:comparison_value_5bus}, the average out-of-sample performances of A4 are the worst among all the aforementioned approaches regardless the choice of $\rho$, as A4 only performs the best when the uncertainties exactly follow the empirical distribution generated by the $N$ samples and any deviation of the true distribution from the empirical distribution will influence its performance significantly. However, the average performances of A2 and A3 are not consistent under different correlation coefficient settings, say A2 outperforms A3 in the low correlation regime and A3 beats A2 in the high correlation scenario. In the high-correlation regime, A3 with moment ambiguity set becomes finding a worst-case distribution among all univariate distributions with given mean and variance. This is a relatively small set of distributions, and considering the underlying distribution is Gaussian, the solution yielding from moment ambiguity set is not overly conservative, and can effectively identify an OPF strategy close to the true optimal one. On the contrary, A2 with the Wasserstein ambiguity set may hedge against some distributions unlikely to happen, which makes the decision over-conservative. In the low-correlation regime, the samples are not so concentrating, so uncertainties towards any direction does not affect the correlation structure too much, hence Wasserstein set along performs well. But moment ambiguity set is too conservative since it contains too many distributions. Thus A2 outperforms A3. As a hybrid approach, A1 takes advantages of A2 and A3, and performs consistently the best among all the approaches regardless the value of $\rho$.
\begin{table*}[ht!]
\footnotesize
  \centering
  \caption{Objective, dipatch cost and out-of-sample performance of A1-A4 for the 5-bus test system.}\label{Tab:comparison_value_5bus}
  \begin{tabular}{ccccccccccc}
  \toprule
  \multirow{2}*{$\rho$}&\multicolumn{2}{c}{\multirow{2}*{}}& \multicolumn{2}{c}{W$\&$M-DROPF} & \multicolumn{2}{c}{W-DROPF} & \multicolumn{2}{c}{M-DROPF} & \multicolumn{2}{c}{SAA}\\
  \cmidrule(lr){4-5}\cmidrule(lr){6-7}\cmidrule(lr){8-9}\cmidrule(lr){10-11}
  & & & DLR & SLR & DLR & SLR & DLR & SLR & DLR & SLR\\
  \midrule
  \multirow{6}*{0.4}&\multirow{3}*{DC ($\$$)} & avg & $5.497\times10^3$ & $6.322\times10^3$ &$5.425\times10^3$ & $6.212\times10^3$ &$5.551\times10^3$ & $6.439\times10^3$ &$5.452\times10^3$ & $6.379\times10^3$\\
  & & max & $6.664\times10^3$& $7.797\times10^3$ & $6.714\times10^3$ & $7.923\times10^3$ &$7.232\times10^3$ & $8.539\times10^3$ & $7.055\times10^3$ & $8.219\times10^3$\\
  & & min & $4.503\times10^3$& $5.043\times10^3$ & $4.557\times10^3$ & $5.195\times10^3$ &$4.696\times10^3$ & $5.494\times10^3$ & $4.653\times10^3$ & $5.384\times10^3$ \\
  \cmidrule(lr){2-11}
  &\multirow{3}*{OP ($\$$)} & avg & $7.336\times10^3$ & $7.703\times10^3$ &$7.408\times10^3$ & $7.815\times10^3$ &$7.491\times10^3$ & $7.661\times10^3$ &$7.775\times10^3$ & $7.802\times10^3$\\
  && max &$7.724\times10^3$ & $8.265\times10^3$  & $7.837\times10^3$& $8.464\times10^3$ & $8.099\times10^3$& $8.504\times10^3$ & $8.372\times10^3$& $8.958\times10^3$ \\
  && min &$7.120\times10^3$ & $6.978\times10^3$  & $7.136\times10^3$& $7.065\times10^3$ & $7.161\times10^3$& $7.304\times10^3$ & $7.251\times10^3$& $7.324\times10^3$ \\
  \midrule
  \multirow{6}*{0.9}&\multirow{3}*{DC ($\$$)} & avg & $5.688\times10^3$ & $6.257\times10^3$ &$6.151\times10^3$ & $6.212\times10^3$ &$5.661\times10^3$ & $6.340\times10^3$ &$5.794\times10^3$ & $6.663\times10^3$\\
  & & max & $6.989\times10^3$& $6.919\times10^3$ & $7.300\times10^3$ & $7.923\times10^3$ &$7.665\times10^3$ & $7.842\times10^3$ & $7.366\times10^3$ & $8.103\times10^3$\\
  & & min & $4.966\times10^3$& $5.661\times10^3$ & $4.241\times10^3$ & $5.195\times10^3$ &$4.580\times10^3$ & $4.762\times10^3$ & $4.154\times10^3$ & $4.362\times10^3$ \\
  \cmidrule(lr){2-11}
  &\multirow{3}*{OP ($\$$)} & avg & $7.253\times10^3$ & $7.471\times10^3$ &$7.423\times10^3$ & $7.720\times10^3$ &$7.341\times10^3$ & $7.708\times10^3$ &$7.505\times10^3$ & $7.543\times10^3$\\
  && max &$7.782\times10^3$ & $7.626\times10^3$  & $8.178\times10^3$& $8.342\times10^3$ & $8.005\times10^3$& $8.245\times10^3$ & $8.478\times10^3$& $8.563\times10^3$ \\
  && min &$7.064\times10^3$ & $7.135\times10^3$  & $7.119\times10^3$& $6.905\times10^3$ & $7.072\times10^3$& $7.213\times10^3$ & $7.159\times10^3$& $6.873\times10^3$ \\
  \bottomrule
  \end{tabular}
\end{table*}

\begin{figure}
\centering
\subfigure[Outputs of generators.]{\label{fig:low_gen}
\includegraphics[width=0.22\textwidth]{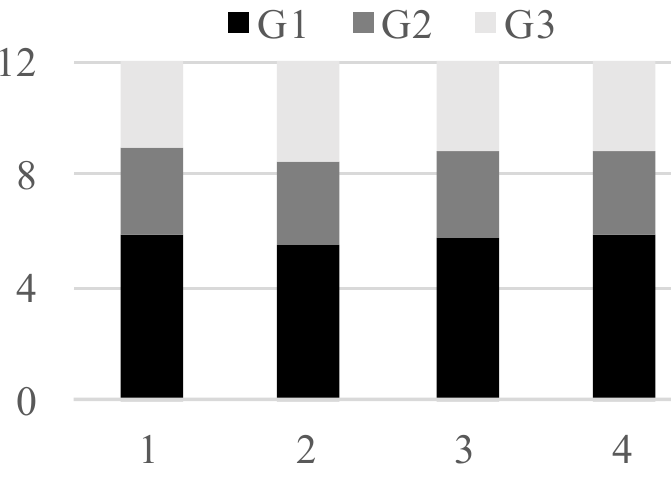}}
\subfigure[Participation factors.]{ \label{fig:low_alpha}
\includegraphics[width=0.22\textwidth]{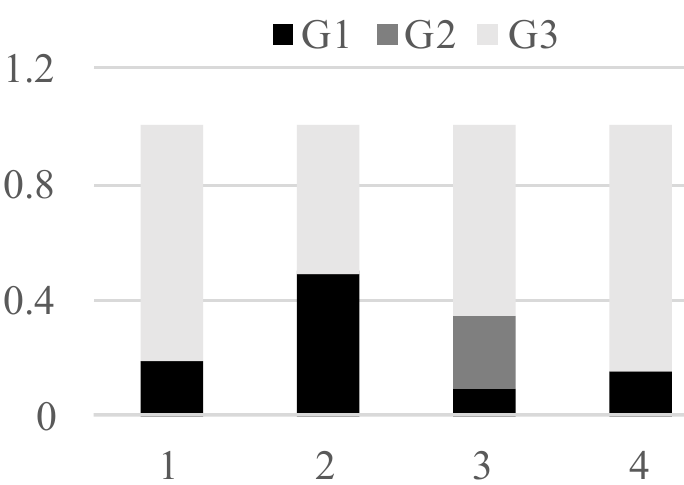}} \\
\subfigure[Upward reserves.]{\label{fig:low_res_pos}
\includegraphics[width=0.22\textwidth]{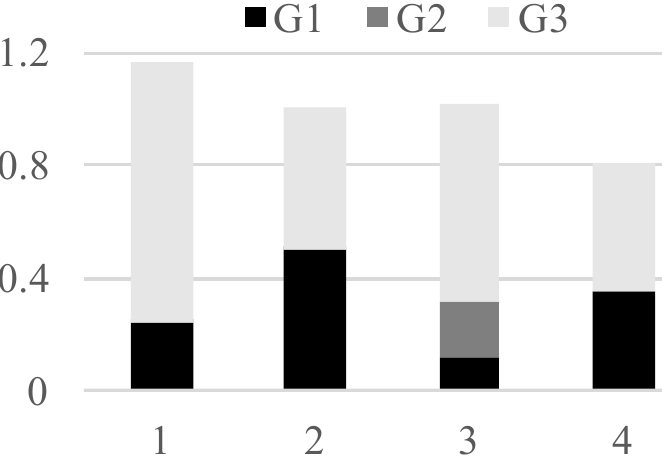}}
\subfigure[Downward reserves.]{ \label{fig:low_res_neg}
\includegraphics[width=0.22\textwidth]{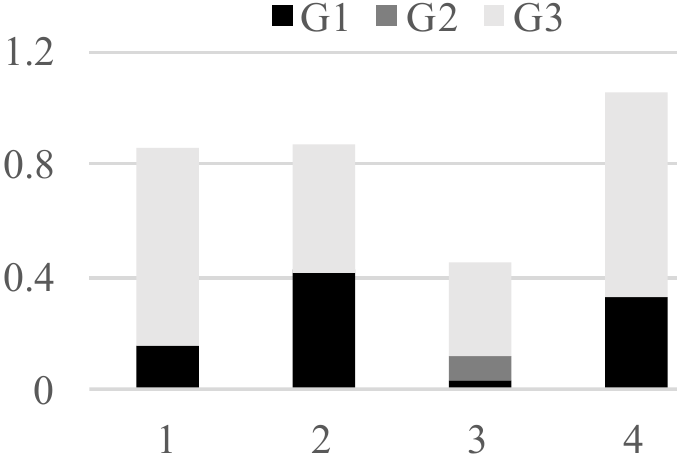}}
\caption{Average OPF strategies of A1-A4 when $\rho=0.4$.} \label{fig:low_results}
\end{figure}

\begin{figure}
\centering
\subfigure[Outputs of generators.]{\label{fig:high_gen}
\includegraphics[width=0.22\textwidth]{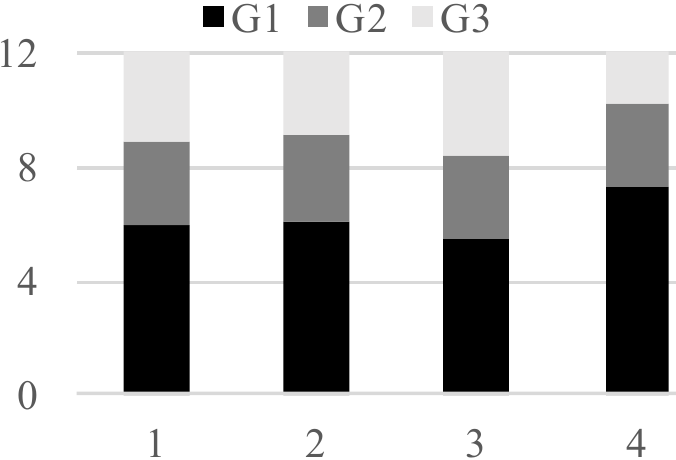}}
\subfigure[Participation factors.]{ \label{fig:high_alpha}
\includegraphics[width=0.22\textwidth]{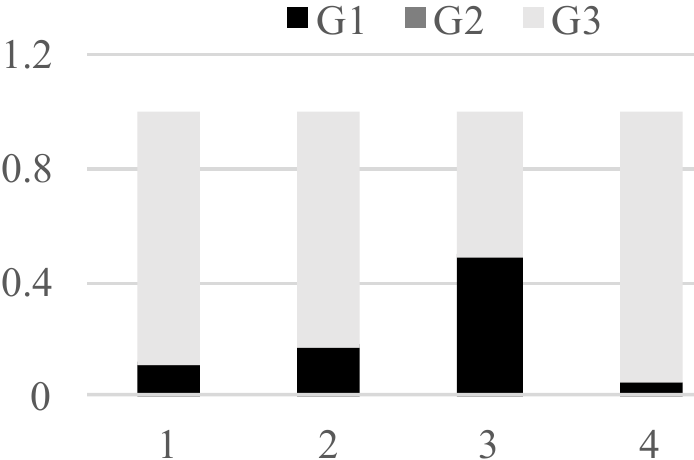}} \\
\subfigure[Upward reserves.]{\label{fig:high_res_pos}
\includegraphics[width=0.22\textwidth]{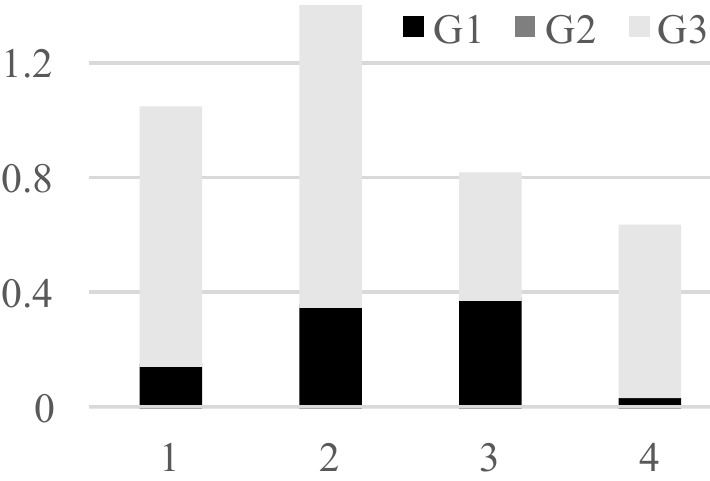}}
\subfigure[Downward reserves.]{ \label{fig:high_res_neg}
\includegraphics[width=0.22\textwidth]{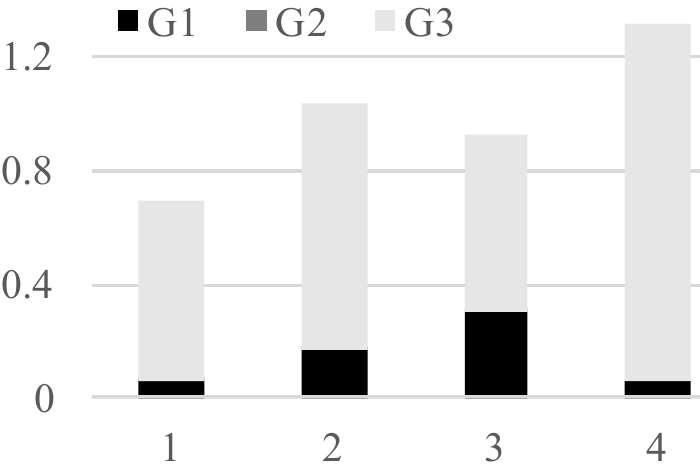}}
\caption{Average OPF strategies of A1-A4 when $\rho=0.9$.} \label{fig:high_results}
\end{figure}

The average OPF strategies of A1-A4 when $\rho$ is set as $0.4$ and $0.9$ are shown in Figs. \ref{fig:low_results} and \ref{fig:high_results}, respectively, where the differences among the performances of the aforementioned approaches can be further explained. Let us first look at the low-correlation case. From Fig. \ref{fig:low_gen}, the outputs of generators are almost the same in A1-A4. However, A1 aims to purchase more upward reserves as the penalty coefficient of load shedding are the highest, which can be viewed in Fig. \ref{fig:low_res_pos} and is the major reason for its best average performances. Both the upward and downward reserves of A2 are less than A1, which accounts for its relatively lower average dispatch costs and slightly poorer performance than A1. A4 purchases the least amount of upward reserves and buys too much downward reserves, resulting in the second highest average dispatch costs and the worst average performances. A3 chooses to commit reserves from $\text{G}_2$, the most expensive unit among all, and does not buy enough downward reserve, which explain its highest average dispatch costs and second worst average performances. Then the high-correlation scenario is visited. From Figs. \ref{fig:high_res_pos} and \ref{fig:high_res_neg}, A2 gets the most upward reserves and the second most downward reserves, which can significantly reduces the operational risk, but the highest dispatch costs partly hinder its average performances. And apparently A4 does not get enough upward reserves, leading to the worst average performances again. Compared with A3, A1 buys more upward reserves to reduce the load shedding risk and commits less downward reserves to sacrifice its performance in wind curtailment, as the former has a larger penalty coefficient in this case. Though the operational risk of A1 might be larger than A2, the OPF strategies are much less conservative, giving rise to its best average performances once again.

\subsection{Effectiveness of Approximation Techniques}
In the sequel, the effectiveness of the DROPF approximation techniques are analyzed. The sample generation and parameter tuning approaches are the same as Section IV.B. The simulation is repeated by 100 times, one half with $\rho=0.4$ and the other half with $\rho=0.9$. The numerical and computational performances of the proposed model (\ref{eq:Risk_term}) and approximated DROPF models (\ref{eq:DROPF_approximate_small}) and (\ref{eq:DROPF_approximate_big}) are compared under 3 different scales of sample dataset, i.e. $N=10,20,50$. The numbers of auxiliary constraints and variables of A1 under the aforementioned settings are presented in Table \ref{Tab:computational_scale_5bus}. And the results are summarized in Table \ref{Tab:approximation_5bus}.

\begin{table}[ht!]
\footnotesize
  \centering
  \caption{Number of auxiliary constraints and variables of A1: the 5-bus system.}\label{Tab:computational_scale_5bus}
  \begin{tabular}{ccccc}
   \toprule
   \multirow{2}*{Sample Number} &\multirow{2}*{Penalties} & \multirow{2}*{Auxiliary Variable} & \multicolumn{2}{c}{Auxiliary Constraint}  \\ 
   \cmidrule(lr){4-5}
   & & & SDP & Linear\\
   \midrule
   \multirow{3}*{$N=10$} & (\ref{eq:Risk_term})  & $69132$ & $17281$ & $17281$ \\
   & (\ref{eq:DROPF_approximate_small})  & $1756$ & $433$ & $433$\\
   & (\ref{eq:DROPF_approximate_big}) & $948$ & $219$ & $219$ \\
   \midrule
   \multirow{3}*{$N=20$}& (\ref{eq:Risk_term}) & $138262$ & $34561$ & $34561$ \\
   & (\ref{eq:DROPF_approximate_small}) & $3506$ & $863$ & $863$ \\
   & (\ref{eq:DROPF_approximate_big}) & $1878$ & $429$ & $429$\\
   \midrule
   \multirow{3}*{$N=50$}& (\ref{eq:Risk_term}) &$345652$& $86401$&$86401$  \\
   & (\ref{eq:DROPF_approximate_small}) & $8756$ & $2153$ & $2153$\\
   & (\ref{eq:DROPF_approximate_big}) & $4668$ & $1059$ & $1059$\\
   \bottomrule
  \end{tabular}
\end{table}

From Table \ref{Tab:approximation_5bus}, it can be observed that the approximated models can always bring significant computational benefits, say the computational time of (\ref{eq:DROPF_approximate_big}) is nearly 1/100 of (\ref{eq:Risk_term}), in all the cases, due to the SDP constraint reduction by approximating the risk-related terms, which can be observed from the last column of Table \ref{Tab:computational_scale_5bus}. As a result, the average performances of approximated models might be sightly worse than accurate one, on account of the additional conservativeness, which can be learnt from the performances of the cases with $N=20, 50$. However, when the samples dataset is small, say $N=10$, the approximated models may outperform the accurate one, as the information revealed from the ambiguity set might be insufficient, leading to an over optimistic OPF strategy. It should be noted that all three groups of samples with $N=10,20,50$ are independently generated, which means the operation costs among different columns of Table \ref{Tab:approximation_5bus} are not comparable.
\begin{table}[ht!]
\footnotesize
  \centering
  \caption{Simulation results of approximated models.}\label{Tab:approximation_5bus}
  \begin{tabular}{ccccc}
   \toprule
   \multirow{2}*{}& \multirow{2}*{Penalties} & \multicolumn{3}{c}{Sample Number}\\
   \cmidrule(lr){3-5}
   &  & 10 & 20 & 50\\
  \midrule
  \multirow{3}*{OP ($\$$)} & (\ref{eq:Risk_term}) & $7.327\times10^3$  &$7.305\times10^3$ & $7.271\times10^3$ \\
  & (\ref{eq:DROPF_approximate_small}) & $7.289\times10^3$&$7.309\times10^3$ &$7.289\times10^3$ \\
  & (\ref{eq:DROPF_approximate_big}) &$7.272\times10^3$ & $7.314\times10^3$ & $7.310\times10^3$\\
  \midrule
  \multirow{3}*{Time ($s$)} & (\ref{eq:Risk_term}) & $6.771$& $12.17$ & $40.11$\\
  & (\ref{eq:DROPF_approximate_small}) & $0.212$& $0.776$ &$1.164$ \\
  & (\ref{eq:DROPF_approximate_big}) & $0.083$ & $0.225$ & $0.524$\\
  \bottomrule
  \end{tabular}
\end{table}

\subsection{Impacts of Penalty Coefficients}
In the sequel, the impacts of penalty coefficients in objective function (\ref{eq:DR_Obj}) on the OPF strategy are analyzed. Table \ref{Tab:penalty_5bus_sensitivity} demonstrates the settings for penalty coefficients. Three cases, denoted as cases 1, 2, 3, respectively, are performed for comparison, where case 1 is the benchmark case and the penalty coefficients are approximately one order of magnitude higher than the cost parameters of generators, and the penalty coefficients of case 2 and case 3 are one order of magnitude higher and lower than case 1, respectively. Then, we repeat the simulations on the proposed W$\&$M-DROPF model with $N=20$, where the samples in all these cases are the same.
\begin{table}[ht!]
\footnotesize
  \centering
  \newcommand{\tabincell}[2]{\begin{tabular}{@{}#1@{}}#2\end{tabular}}
  \caption{Penalty Coefficients: the 5-bus Test System.}\label{Tab:penalty_5bus_sensitivity}
  \begin{tabular}{cccc}
  \toprule
   &$\beta_d$ ($\$$/MWh) & $\beta_w$ ($\$$/MWh) & $\beta_l$ ($\$$/MWh)\\
   \midrule
   Case 1&$1\times 10^4$ & $1\times 10^3$ & $5\times 10^3$\\
   Case 2&$1\times 10^5$ & $1\times 10^4$ & $5\times 10^4$\\
   Case 3&$1\times 10^3$ & $1\times 10^2$ & $5\times 10^2$\\
  \bottomrule
  \end{tabular}
\end{table}
\begin{table}[ht!]
\footnotesize
  \centering
  \newcommand{\tabincell}[2]{\begin{tabular}{@{}#1@{}}#2\end{tabular}}
  \caption{Results under different penalty coefficients.}\label{Tab:penalty_5bus_results}
  \begin{tabular}{ccccccc}
  \toprule
   & DC ($\$$) & & \tabincell{c}{Output\\ (MW)}& \tabincell{c}{UR\\ (MW)}& \tabincell{c}{DR \\(MW)}\\
   \midrule
   \multirow{3}*{Case 1} & \multirow{3}*{$6.171\times 10^3$} & $G_1$ & $7.41$ & $0.10$ & $0.06$\\
   & & $G_2$ & $2.24$ & $0.26$ & $0.11$ \\
   & & $G_3$ & $2.36$ & $1.07$ & $0.56$ \\
   \midrule
   \multirow{3}*{Case 2} & \multirow{3}*{$7.4462\times 10^3$} & $G_1$ & $7.86$ & $0.94$ & $0.88$\\
   & & $G_2$ & $2.34$ & $0.66$ & $0.75$ \\
   & & $G_3$ & $1.80$ & $0$ & $0$ \\
   \midrule
   \multirow{3}*{Case 3} & \multirow{3}*{$3.9113\times 10^3$} & $G_1$ & $7.16$ & $0$ & $0$\\
   & & $G_2$ & $2.18$ & $0$ & $0$ \\
   & & $G_3$ & $2.66$ & $0$ & $0$ \\
  \bottomrule
  \end{tabular}
\end{table}

The simulation results are summarized in Table \ref{Tab:penalty_5bus_results}, where DC in the second column is short for dispatch costs (summation of generation costs and reserve costs), and UR and DR are short for upward reserve and downward reserve, respectively. It can be observed that dispatch costs have a positive relationship with the penalty coefficients, say the larger penalty coefficients are, the large dispatch costs would be. Similarly, from Table \ref{Tab:penalty_5bus_results}, the penalty coefficients also have significant impacts on the outputs and allocated reserves of generators. Specifically, in case 3, the penalty coefficients are about the same with cost parameters of generators, no reserve would be committed. And in case 2, the penalty coefficients are far greater (two orders of magnitude higher) than the generator cost parameters, both the total upward and downward reserves are in relatively high levels. Meanwhile, no reserve would be committed from $G_3$ though it has the lowest prices, as its location is relatively far from the demands and more lines will be used to deliver the reserves, indicating additional line overload risk if reserves are committed from it.

\subsection{Impacts of Sample Number}
Though the number of available historical data to construct the empirical distribution $\nu$ is assumed to be relatively small in our case, which has been mentioned in Remark 3 of Section II.D and is also the major motivation of proposing the DROPF models, the performances of A1-A4 with larger sample sets, say $N=100$ and $N=1000$ are demonstrated in this sequel. The sample generation and parameter tuning approach are identical to the previous subsections. The simulation is repeated by 100 times, where 50 times with $\rho=0.4$ and the other 50 times with $\rho=0.9$. The penalty terms in A1-A3 are approximated by (\ref{eq:DROPF_approximate_big}) to reduce the additional computational burden brought by the relatively large numbers of samples. The average out-of-sample performances of A1-A4 are gathered in Table \ref{Tab:results_5bus_N100}.
\begin{table}[ht!]
\footnotesize
  \centering
  \caption{Average performances of A1-A4 in the 5-bus test system with larger sample sets.}\label{Tab:results_5bus_N100}
  \begin{tabular}{ccccc}
   \toprule
   \multirow{2}*{$N$} & \multicolumn{4}{c}{Out-of-sample Performance ($\$$)}\\
   \cmidrule(lr){2-5}
   & A1 & A2 & A3 & A4 \\
   \midrule
   $100$ & $7.056\times 10^3$ & $7.093\times 10^3$ & $7.127\times 10^3$ & $7.136\times 10^3$\\
   $1000$& $6.920\times 10^3$ & $6.928\times 10^3$ & $6.939\times 10^3$ & $6.941\times 10^3$\\
   \bottomrule
  \end{tabular}
\end{table}

From Table \ref{Tab:results_5bus_N100}, it can be observed that A1 still performs the best among all the listed OPF models in both $N=100$ and $N=1000$ cases, as it combines the merits of A2 and A3. Moreover, the performances of A3 become the second worst, as only limited information, say the mean value and covariance matrix, of the samples are used. Though A4 still has the worst out-of-sample performances in both cases, the relative gap between the performances of A1 and A4 decreases from $1.13\%$ to $0.30\%$ when $N$ increases from $100$ to $1000$, which indicates SAA would get similar performances with the listed DROPF models when the sample set is sufficient large.

\section{Case Study}
In this section, we present numerical experiments on two larger test systems, which are the IEEE 118-bus system and the Polish 2736-bus power system, to demonstrate the performance of the proposed model and algorithm. The simulation environment and solver settings are identical to Section IV.

\subsection{Simulation Results of the IEEE 118-bus System}
To demonstrate the scalability and efficiency of the proposed methods, they are applied to a larger system, consisting of a modified IEEE 118-bus test system. The test system has 54 generators and 186 transmission lines. Three wind farms are connected to the system at buses 17, 66 and 94. Please refer to \cite{5bus} for the topology and system data. For the ease of analysis, only 30 generators participate in AGC, and 59 lines are left after inactive transmission capacity constraints identification. The simulation is repeated by 100 times with $N=10,20,50$, where the risk-related terms in A1-A3 are approximated by (\ref{eq:DROPF_approximate_big}). The major computational burden of A1-A4 are listed in Table \ref{Tab:computational_scale_118bus}. From Table \ref{Tab:computational_scale_118bus}, the computational burden of A1 is the highest due to the SDP constraints, whose number grows approximately linearly with respect to the sample number, and the computational burden of A3 does not change with the sample number, as only the mean and variance of the samples are needed. Besides, the computational burden of A2 and A4 are almost the same, and both of them are QPs.
\begin{table}[ht!]
\footnotesize
  \centering
    \caption{Number of auxiliary constraints and variables of A1-A4: the IEEE 118-bus system.}\label{Tab:computational_scale_118bus}
  \begin{tabular}{ccccc}
   \toprule
   \multicolumn{2}{c}{\multirow{2}*{Model Type}} & \multirow{2}*{Auxiliary Variable} & \multicolumn{2}{c}{Auxiliary Constraint}  \\ 
   \cmidrule(lr){4-5}
   & & & SDP & Linear\\
   \midrule
   \multirow{3}*{A1: W$\&$M-DROPF} & $N=10$  & $185568$ & $3089$ & $3089$ \\
   & $N=20$  & $370898$ & $6059$ & $6059$\\
   & $N=50$ & $926888$ & $14969$ & $14969$ \\
   \midrule
   \multirow{3}*{A2: W-DROPF}& $N=10$ & $1396$ & $0$ & $3148$ \\
   & $N=20$ & $2199$ & $0$ & $6118$ \\
   & $N=50$ & $6069$ & $0$ & $15028$\\
   \midrule
   \multicolumn{2}{c}{A3: M-DROPF} & $464933$& $416$ & $0$ \\
   \midrule
   \multirow{3}*{A4: SAA}& $N=10$ &$1190$& $0$&$2970$  \\
   & $N=20$ & $2380$ & $0$ & $5940$\\
   & $N=50$ & $5950$ & $0$ & $14850$\\
   \bottomrule 
  \end{tabular}
\end{table}

The average performances of A1-A4 are shown in Table \ref{Tab:results_118bus}. From Table \ref{Tab:results_118bus}, it can be observed that A1 still performs the best among all the approaches in all the cases in terms of average performances, even though its risk-related terms are approximated and additional conservativeness might be introduced. The performances of A2 and A3 are almost the same considering the mixture of high- and low-correlation sample datasets. Still, A4 has the worst average performances. For the computational time, A1 is the most time consuming, as a large number of SDPs are tackled. However, the computational time of A1 is still acceptable for a moderate-size power system, considering the performance advantage over the other approaches. Likewise, there is no comparability among the operation costs of one model with different sample numbers. However, it can be observed that relative gap between the average performances of A1 and A4 decreases when sample number grows, which is in consistent with the observation in Section IV.F, indicating the average performances of SAA and the proposed data-driven approaches might be close when sample number is sufficient large.
\begin{table}[ht!]
\footnotesize
  \centering
  \caption{Simulation results of the IEEE 118-bus test system.}\label{Tab:results_118bus}
  \begin{tabular}{cccc}
  \toprule
  Model Type & Sample Number & OP ($\$$) & Time ($s$)\\
  \midrule
  \multirow{3}*{A1} & $N=10$ & $1.235\times10^5$ & $1.873$\\
  & $N=20$ & $1.187\times10^5$ & $4.227$\\
  & $N=50$ & $1.226\times10^5$ & $14.77$\\
  \midrule
  \multirow{3}*{A2} & $N=10$ & $1.269\times10^5$ & $0.052$\\
  & $N=20$ & $1.215\times10^5$ & $0.227$\\
  & $N=50$ & $1.249\times10^5$ & $0.938$\\
  \midrule 
  \multirow{3}*{A3} & $N=10$ & $1.277\times10^5$ & $0.463$\\
  & $N=20$ & $1.213\times10^5$ & $0.511$\\
  & $N=50$ & $1.255\times10^5$ & $0.489$\\
  \midrule
  \multirow{3}*{A4} & $N=10$ & $1.337\times10^5$ & $0.022$\\
  & $N=20$ & $1.251\times10^5$ & $0.211$\\
  & $N=50$ & $1.263\times10^5$ & $0.852$\\
  \bottomrule
  \end{tabular}
\end{table}

\subsection{Simulation Results of the Polish 2736-bus Power System}
In this sequel, we test the proposed methods on the Polish 2736-bus power system in the Matpower toolbox, which is the Polish 400, 220 and 110 kV networks during summer 2004 peak conditions. The test system has 420 generators and 3504 transmission lines. Please refer to \cite{Polish} for the topology and system data. Ten 400MW wind farms are connected to the system at buses 1 to 10. The generation and reserve costs parameters are assigned with the ones of the IEEE 118-bus system. We assume all the generators participate in AGC and all the transmission lines are monitored by DLR devices. Similarly, three groups of samples are generated with $N=10,20,50$. Before demonstrating the simulation results, the numbers of auxiliary variables and constraints generated by A1-A4 for the Polish power system are listed as below
\begin{table}[ht!]
\footnotesize
  \centering
    \caption{Number of auxiliary constraints and variables of A1-A4: the Polish 2736-bus power system.}\label{Tab:computational_scale_Polish}
  \begin{tabular}{ccccc}
   \toprule
   \multicolumn{2}{c}{\multirow{2}*{Model Type}} & \multirow{2}*{Auxiliary Variable} & \multicolumn{2}{c}{Auxiliary Constraint}  \\ 
   \cmidrule(lr){4-5}
   & & & SDP & Linear\\
   \midrule
   \multirow{3}*{A1} & $N=10$  & $4.285\times 10^8$ & $126264$ & $126264$ \\
   & $N=20$  & $8.569\times 10^8$ & $248184$ & $248184$\\
   & $N=50$ & $2.142\times 10^9$ & $613944$ & $613944$ \\
   \midrule
   \multirow{3}*{A2}& $N=10$ & $47784$ & $0$ & $129768$ \\
   & $N=20$ & $91224$ & $0$ & $251688$ \\
   & $N=50$ & $221544$ & $0$ & $617448$\\
   \midrule
   \multicolumn{2}{c}{A3} & $5.366\times 10^{10}$& $16536$ & $0$ \\
   \midrule
   \multirow{3}*{A4}& $N=10$ &$43440$& $0$&$121920$  \\
   & $N=20$ & $86880$ & $0$ & $243840$\\
   & $N=50$ & $217200$ & $0$ & $609600$\\
   \bottomrule 
  \end{tabular}
\end{table}
 
From Table \ref{Tab:computational_scale_Polish}, it can be inferred that both A1 and A3 might not be handled by the current computation platform, which is a laptop with a 2.2 GHz CPU and a 4GB memory, as their numbers of auxiliary variables are relatively large and both of them require solving SDPs. Then the simulation is repeated by 10 times and the results are summarized in Table \ref{Tab:results_Polish}.
\begin{table}[ht!]
\footnotesize
  \centering
  \caption{Simulation results of the Polish 2736-bus power system.}\label{Tab:results_Polish}
  \begin{tabular}{cccc}
  \toprule
  \multirow{2}*{Sample Number}& &\multicolumn{2}{c}{Model Type}\\
  \cmidrule(lr){3-4}
  & & A2 & A4\\
  \midrule
  \multirow{2}*{$N=10$}&OP ($\$$) & $2.327\times10^7$& $2.447\times10^7$\\
  &Time ($s$)  & $322$& $285$\\
  \midrule
  \multirow{2}*{$N=20$}&OP ($\$$) & $2.439\times10^7$& $2.501\times10^7$\\
  &Time ($s$)  & $419$& $347$\\
  \midrule
  \multirow{2}*{$N=50$}&OP ($\$$) & $2.455\times10^7$& $2.493\times10^7$\\
  &Time ($s$) & $761$& $646$\\
  \bottomrule
  \end{tabular}
\end{table} 

It should be noted that the performances of A1 and A3 are not shown in Table \ref{Tab:results_Polish}. The reason is that the simulation software, which is Matlab in our case, would encounter out of memory issue when A1 or A3 is performed in all three sample groups, indicating the intractability of the proposed DROPF or the first two orders moment based DROPF models for practical large-scale system on personal computers. However, both A2 and A4 can be solved within an acceptable time, which can be observed from Table \ref{Tab:results_Polish}. Meanwhile, the average performances of A2 are always better than A4 in all three cases, validating the effectiveness of the proposed W-DROPF model.

\section{Conclusion}
DLR has been proved to be of great value to maximize the capability of power system operation to hedge against uncertainties of outputs of wind generation, facilitating the utilization of wind generation simultaneously. However, the implementation of DLR will introduce additional uncertainties, as the actual line rating cannot be accurately known or predicted beforehand. To address this issue, a risk-based DROPF model with DLR is proposed. Both the second-order moment ambiguity set and Wasserstein ambiguity set are considered in the proposed model to better capture the correlation of samples and preserve the robustness of OPF strategy to rare samples. The original min-max W$\&$M-DROPF model is reformulated as a single level minimization problem with exponential many SDP constraints based on strong duality theory, which is readily to solve by the off-the-shelf solvers. A mild approximation of risk terms in the proposed model is derived to reduce the computational burden. For practical large-scale test systems, the underlying computational burden might no be well handled by small computation platforms, such as personal computers, even if the risk-term approximation is adopted. The Wasserstein distance constrained DROPF model is prepared for this situation, whose tractable reformulation is a QP mathematically.

Simulation results show that the proposed W$\&$M-DROPF model has better out-of-sample performances over M-DROPF model (the one with second-order moment constrained ambiguity set), W-DROPF model (the one with Wasserstein distance constrained ambiguity set), and of course the SAA approach. And the advantage still exists even if its approximation is dealt with rather than the original one, which is also revealed by the simulation results. For practical large-scale test systems, such as the Polish 2736-bus power system, though the W$\&$M-DROPF and M-DROPF models cannot be handled due to the heavy computational burden of their reformulations, W-DROPF still outperforms SAA in terms of out-of-sample performance, which offers an effective and efficient alternative. In fact, the ambiguity set constructed in this paper can be applied in many other power system decision-making problems, which need to be formulated in a distributionally robust manner, such as economic dispatch and reserve procurement.

\section*{Appendix}
\subsection{Detailed Expressions of Coefficients in (\ref{eq:strong_dual})}
The index set $\{1,\ldots,K\}$ is reparameterized as
\[\begin{aligned}
  &\Big\{(j_{g_d},j_{g_w},j_l):j_{g_d},j_{g_w}\in \{1,2\}, j_\ell\in\{1,2,3\}, \\ & \hspace{12em} \forall g_d,g_w\in \mathcal{G},l\in \mathcal{L}\Big\},
  \end{aligned}
\]
and thus $K=4^G\times 3^L$.
We express
\[
  \vec{a}_{(j_{g_d},j_{g_w},j_l)}(\vec{x}) ~=~ \vec{a}^{g_d}_{j_{g_d}}(\vec{x}) + \vec{a}^{g_w}_{j_{g_w}}(\vec{x}) + \vec{a}^{l}_{j_l}(\vec{x}),
\]
\[
  b_{(j_{g_d},j_{g_w},j_l)}(\vec{x}) ~=~ b^{g_d}_{j_{g_d}}(\vec{x}) + b^{g_w}_{j_{g_w}}(\vec{x}) + b^{l}_{j_l}(\vec{x}),
\]
where $\vec{a}^{g_d}_{j_{g_d}}(\vec{x})$, $\vec{a}^{g_w}_{j_{g_w}}(\vec{x})$, $\vec{a}^{l}_{j_l}(\vec{x})$, $b^{g_d}_{j_{g_d}}(\vec{x})$, $b^{g_w}_{j_{g_w}}(\vec{x})$, $b^{l}_{j_l}(\vec{x})$ are defined through
\[
  \vec{a}^{g_d}_{j_{g_d}}(\vec{x})^\top \vec{\xi} = \begin{cases}
    -\beta_d\alpha_g\sum_{w\in \mathcal{W}}\tilde{p}_w,& \ j_{g_d}=1,  \\
    0,& \ j_{g_d}=2,
  \end{cases}
\]
\[
  \vec{a}^{g_w}_{j_{g_w}}(\vec{x})^\top \vec{\xi} = \begin{cases}
    \beta_w\alpha_g\sum_{w\in \mathcal{W}}\tilde{p}_w,& \ j_{g_w}=1, \\
    0,& \ j_{g_w}=2,
  \end{cases}
\]
\[\begin{aligned}
  & \vec{a}^{l}_{j_l}(\vec{x})^\top \vec{\xi}\\
   =& \begin{cases}
    \beta_l\Big(\sum\limits_{w\in \mathcal{W}}\pi_{wl}\tilde{p}_w-\sum\limits_{g\in \mathcal{G}}\pi_{gl}\alpha_g\sum\limits_{w\in \mathcal{W}}\tilde{p}_w-\tilde{p}_l\Big),& \ j_l=1, \\
    \beta_l\Big(\tilde{p}_l-\sum\limits_{w\in \mathcal{W}}\pi_{wl}\tilde{p}_w+\sum\limits_{g\in \mathcal{G}}\pi_{gl}\alpha_g\sum\limits_{w\in \mathcal{W}}\tilde{p}_w\Big),& \ j_l=2, \\
    0, & \ j_l=3,
  \end{cases}
  \end{aligned}
\]
and
\[
  b^{g_d}_{j_{g_d}}(\vec{x}) = \begin{cases}
    -\beta_dr_g^++\beta_d\alpha_g\sum_{w\in \mathcal{W}}p_w,& \ j_{g_d}=1,  \\
    0,& \ j_{g_d}=2,
  \end{cases}
\]
\[
  b^{g_w}_{j_{g_w}}(\vec{x}) = \begin{cases}
    -\beta_wr_g^--\beta_w\alpha_g\sum_{w\in \mathcal{W}}p_w,& \ j_{g_w}=1, \\
    0,& \ j_{g_w}=2,
  \end{cases}
\]
\[\begin{aligned}
  & b^{l}_{j_l}(\vec{x})\\
  =& \begin{cases}
    \displaystyle \beta_l\sum_{g\in \mathcal{G}}\pi_{gl}\Big(p_g+\alpha_g\sum_{w\in \mathcal{W}}p_w\Big)-\beta_l\sum_{d\in \mathcal{D}}\pi_{dl}p_d,& j_l=1, \\
    \displaystyle-\beta_l\sum_{g\in \mathcal{G}}\pi_{gl}\Big(p_g+\alpha_g\sum_{w\in \mathcal{W}}p_w\Big)+\beta_l\sum_{d\in \mathcal{D}}\pi_{dl}p_d,& j_l=2, \\
    0, & j_l=3.
  \end{cases}
  \end{aligned}
\]

\subsection{Tractable Reformulations of M-DROPF}
The M-DROPF model has a strong dual reformulation as follows (Theorem 4 in \cite{delage2010distributionally}).
\begin{subequations}\label{eq:DROPF_only_moment}
\begin{align}
& \min_{\substack{\vec{x},\lambda\\ \vec{\Gamma}\succeq0,\vec{\zeta}}}\ \vec{f}(\vec{x})+\tau\text{tr}\left(\vec{\Gamma}\hat{\vec{\Sigma}}\right)+\hat{\vec{m}}^\top\vec{\Gamma}\hat{\vec{m}}+\hat{\vec{m}}^\top\vec{\zeta}+\lambda\label{eq:DROPF_moment_obj}\\
\ s.t.\ \ & \vec{Ax}\le\vec{h},\label{eq:moment_linear}\\
\begin{split}
&\left[\begin{array}{cc}
            \vec{\Gamma} & -\vec{a}_k(\vec{x})/2+\vec{\zeta}/2\\
            (-\vec{a}_k(\vec{x})/2+\vec{\zeta}/2)^\top & -b_k(\vec{x})+\lambda
          \end{array}\right] \succeq 0,\\
&\forall 1\leq k\leq K,
\end{split}\label{eq:moment_SDP}
\end{align}
\end{subequations}
where $\vec{x}=[(p_g;r_g^+;r_g^-;\alpha_g): g\in \mathcal{G}]$ is a vector of decision variables; $\lambda$ is an auxiliary variable; $\vec{\zeta}\in\mathbb{R}^{(L+W)}$ is an auxiliary variable vector; $\vec{\Gamma}\in\mathbb{R}^{(L+W)}\times\mathbb{R}^{(L+W)}$ is an auxiliary variable matrix; tr is the trace operator; $\vec{a}_k(\vec{x}),b_k(\vec{x})$ are the coefficients and their detailed expressions can be found in Appendix.A. In (\ref{eq:DROPF_only_moment}), objective function (\ref{eq:DROPF_moment_obj}) is quadratic and convex, constraints (\ref{eq:DROPF_moment_obj}) are SDP constraints, and constraints (\ref{eq:moment_linear}) are linear constraints and are identical to (\ref{eq:DRcom_Con}). Similarly, (\ref{eq:DROPF_only_moment}) is an SDP and is readily to be solved by MOSEK.

\ifCLASSOPTIONcaptionsoff
  \newpage
\fi

\bibliographystyle{IEEEtran}
\bibliography{IEEEabrv,refs_DROPF}
\end{document}